\documentclass[11pt]{amsart}

\usepackage{amsmath,amssymb,amscd}

\usepackage{amsrefs}
\usepackage{enumerate}
\usepackage{xcolor}

\usepackage{amsmath,amssymb}

\usepackage{amsxtra, amsmath}
\usepackage{amssymb, amscd}
\usepackage{graphicx}
\usepackage{mathrsfs}

\usepackage{mathtools}
\usepackage{float}

\newtheorem{thm}{Theorem}[section]
\newtheorem{prop}[thm]{Proposition}

\newtheorem{cor}[thm]{Corollary}

\theoremstyle{definition}
\newtheorem{definition}[thm]{Definition}

\theoremstyle{remark}
\newtheorem{remark}[thm]{Remark}
\newtheorem{example}[thm]{Example}

\usepackage{dsfont}
\allowdisplaybreaks

\title{ Structure of operator algebras for matrix orthogonal polynomials}

\author{Ignacio Bono Parisi}
\author{Ines Pacharoni}



\address{CIEM-FaMAF\\ Universidad Nacional de C\'or\-do\-ba\\
CP 5000, C\'or\-do\-ba,  Argentina}
\email{ignacio.bono@unc.edu.ar, ines.pacharoni@unc.edu.ar}

\begin{document}
\begin{abstract}
In this paper, we study the structure of the differential operator algebra \( \mathcal{D}(W) \) and its associated eigenvalue algebra \( \Lambda(W) \) for matrix-valued orthogonal polynomials. While \( \Lambda(W) \) is isomorphic to \( \mathcal{D}(W) \), its simpler framework allows us to efficiently derive strong results about \( \mathcal{D}(W) \) and its center \( \mathcal{Z}(W) \). We analyze the behavior of the center under Darboux transformations, establishing explicit relationships between the centers of Darboux-equivalent weights. These results are illustrated through the study of both reducible and irreducible matrix weights, including a detailed analysis of an irreducible Jacobi-type weight.

\end{abstract}

\maketitle

\section{Introduction}

The classical Bochner Problem determines all scalar weights \( w(x) \) such that their sequences of orthogonal polynomials are eigenfunctions of a second-order differential operator. The solutions, namely the well-known Hermite, Laguerre, and Jacobi weights, are of fundamental importance in mathematics and physics. These classical orthogonal polynomials play a central role in spectral theory, special functions, and quantum mechanics, serving as canonical examples with remarkable algebraic and analytic properties.

In the matrix-valued setting, the analogous question is to determine matrix weights \( W(x) \) of size \( N \times N \) such that the associated sequence of orthogonal matrix polynomials are eigenfunctions of a second-order matrix differential operator. This problem, often referred to as the {\em Matrix Bochner Problem}, is significantly more intricate due to the richer algebraic structure underlying matrix-valued orthogonal polynomials. 
The problem is equivalent to the existence of a second-order differential operator in the algebra \( \mathcal{D}(W) \), which consists of all differential operators that have the orthogonal polynomials as eigenfunctions.

The first attempt to analyse the whole algebra $\mathcal D(W)$, beyond merely proving the existence of a second-order operator, was carried out in \cite{CG06} for a few concrete weights, with the aid of symbolic computation. A truly structural perspective followed in the pioneering works of Tirao \cite{T11} and Zurrián \cite{Z17}, which provided the first comprehensive descriptions of $\mathcal D(W)$ for explicit examples.

Building on these ideas, 
Casper and Yakimov \cite{CY18} proved that, under certain natural assumptions, every solution to the Matrix Bochner problem can be obtained as Darboux transformations of direct sums of classical scalar weights. Their results established a deep connection between Darboux transformations and the algebra \(\mathcal{D}(W)\), opening the door to new algebraic techniques for studying matrix-valued orthogonal polynomials. Understanding the algebra \( \mathcal{D}(W) \) is therefore fundamental to solving the matrix Bochner problem.

In recent works \cites{BP23, BP24-2}, we constructed explicit families of matrix weights that cannot be obtained as Darboux transformations of classical scalar weights. These examples highlight the need for a more comprehensive understanding of Darboux transformations and their effect on the algebra \( \mathcal{D}(W) \). In particular, Darboux transformations have proven to be a powerful tool for relating the structure of \( \mathcal{D}(W) \) to that of Darboux-equivalent weights \cites{BP23b, BPstrong, BPZ24}.

In this paper, we deepen this analysis by focusing on two key algebraic structures associated with a matrix weight \( W \): the eigenvalue algebra \( \Lambda(W) \) and the center \( \mathcal{Z}(W) \) of \( \mathcal{D}(W) \). The eigenvalue algebra \( \Lambda(W) \), which is isomorphic to \( \mathcal{D}(W) \), provides a simpler framework for working with differential operators. This simplification allows us to derive strong results about the structure of \( \mathcal{D}(W) \) and its center. By leveraging the connection between \( \mathcal{D}(W) \) and \( \Lambda(W) \), we study how Darboux transformations affect \( \mathcal{Z}(W) \) and establish explicit relationships between the centers of Darboux-equivalent weights.

This work reveals the interplay between Darboux transformations, the algebra \( \mathcal{D}(W) \), and the eigenvalue algebra \( \Lambda(W) \). The results provide new insights into the structure of matrix-valued orthogonal polynomials and their associated differential operators, building a bridge between the classical scalar theory and the more complex matrix setting.

\smallskip

The paper is organized as follows: 
In Section \ref{sec-MOP}, we provide the necessary background on matrix-valued orthogonal polynomials and the algebra of differential operators $D(W)$ associated with a matrix weight $W$. 
Section \ref{Darboux-sect} focuses on Darboux transformations for matrix weights. We review their definition and key results, highlighting how Darboux transformations connect the algebras $D(W)$ and $D(\widetilde{W})$ for Darboux-equivalent weights $W$ and $\widetilde{W}$.

In Section \ref{sec-centers}, we study the center $Z(W)$ of the algebra $D(W)$ and analyze its behavior under Darboux transformations. We establish explicit relationships between the centers of Darboux-equivalent weights. 
We prove that weights within the same Darboux-equivalence class of a direct sum of classical weights can be obtained by factoring an operator from the center. This fundamental connection ensures that the centers of their algebras are intrinsically related, establishing a deep interplay between the Darboux transformation and the algebraic structure of the centers.

Section \ref{sect-eigenvalues} introduces the eigenvalue algebra \( \Lambda(W) \), a representation of \( \mathcal{D}(W) \) as matrices with polynomial coefficients in \( \mathbb{C}[n] \). Although \( \Lambda(W) \) and \( \mathcal{D}(W) \) are isomorphic and share the same algebraic structure and properties, working with \( \Lambda(W) \) often proves more convenient. This simpler framework allows us to efficiently derive strong results about the algebra \( \mathcal{D}(W) \), particularly regarding its center \( \mathcal{Z}(W) \). 
We also study the modules $D(W, \widetilde{W})$ of differential operators mapping between polynomials associated with different weights, connecting these results with the structure of $\Lambda(W)$.

In Section \ref{sect-escalar}, we analyze the case where the weight $W$ is a direct sum of classical scalar weights. We explicitly describe the algebra $D(W)$, its center $Z(W)$, and its eigenvalue algebra $\Lambda(W)$, illustrating the non-trivial structure that emerges even in reducible cases.

Finally, in Section \ref{sect-Jacobi}, we study an irreducible Jacobi-type matrix weight \( W \), which arises as a Darboux transformation of a classical diagonal weight. Building upon the results developed in the preceding sections, we determine the algebra \( \mathcal{D}(W) \), its center \( \mathcal{Z}(W) \), and its eigenvalue algebra \( \Lambda(W) \). We provide a set of generators for \( \mathcal{D}(W) \) and describe its structure in detail, 
highlighting the utility and depth of the techniques developed throughout this work.

\

\section{Matrix valued orthogonal polynomials and the algebra $\mathcal D(W)$}\label{sec-MOP}

 Let $W=W(x)$ be a weight matrix of size $N$ on the real line, that is, a complex $N\times N$ matrix-valued smooth function on the interval $(x_0,x_1)$ such that $W(x)$ is positive definite almost everywhere and with finite moments of all orders. Let $\operatorname{Mat}_N(\mathbb{C})$ be the algebra of all $N\times N$ complex matrices and let $\operatorname{Mat}_N(\mathbb{C}[x])$ be the algebra of polynomials in the indeterminate $x$ with coefficients in $\operatorname{Mat}_N(\mathbb{C})$. We consider the following Hermitian sesquilinear form in the linear space $\operatorname{Mat}_N(\mathbb{C}[x])$
\begin{equation*}
  \langle P,Q \rangle =  \langle P,Q \rangle_W = \int_{x_0}^{x_1} P(x) W(x) Q(x)^*\,dx.
\end{equation*}

Given a weight matrix $W$ one can construct sequences 
$\{Q_n\}_{n\in\mathbb{N}_0}$ of matrix-valued orthogonal polynomials, i.e. the 
$Q_n$ are polynomials of degree $n$ with nonsingular leading coefficient and $\langle Q_n,Q_m\rangle=0$ for $n\neq m$.
We observe that there exists a unique sequence of monic orthogonal polynomials $\{P_n\}_{n\in\mathbb{N}_0}$ in $\operatorname{Mat}_N(\mathbb{C}[x])$.
By following a standard argument (see \cite{K49} or \cite{K71}) one shows that the monic orthogonal polynomials $\{P_n\}_{n\in\mathbb{N}_0}$ satisfy a three-term recursion relation
\begin{equation}\label{ttrr}
    x P_n(x)=P_{n+1}(x) + B_{n}P_{n}(x)+ C_nP_{n-1}(x), \qquad n\in\mathbb{N}_0,
\end{equation}
where $P_{-1}=0$ and $B_n, C_n$ are matrices depending on $n$ and not on $x$.

The introduction of matrix-valued weights naturally leads to the study of differential operators that act on these matrix-valued orthogonal polynomials. While in the scalar case, the algebra of differential operators is typically generated by a single second-order operator, the matrix-valued case exhibits a much richer structure. 

Throughout this paper, we consider that an arbitrary matrix differential operator
\begin{equation}\label{D2}
  {D}=\sum_{i=0}^m \partial ^i F_i(x),\qquad \partial=\frac{d}{dx},
\end{equation}
acts on the right on a matrix-valued function $P$ i.e.
$(P\cdot D)(x)=\sum_{i=0}^m \partial ^i (P)(x)F_i(x).$

We consider the algebra of these operators with polynomial coefficients $$\operatorname{Mat}_{N}(\Omega[x])=\Big\{D = \sum_{j=0}^{n} \partial^{j}F_{j}(x) \, : F_{j} \in \operatorname{Mat}_{N}(\mathbb{C}[x]) \Big \}.$$
\noindent 
More generally, when necessary, we will also consider $\operatorname{Mat}_{N}(\Omega[[x]])$, the set of all differential operators with coefficients in  $\mathbb{C}[[x]]$, the ring of power series with coefficients in $\mathbb{C}$.

\

Given a weight matrix $W$, we can consider  the associative algebra
\begin{equation}\label{algDW}
  \mathcal D(W)=\Big\{D\in \operatorname{Mat}_{N}(\Omega[x]) :  \,P_n \cdot D=\Lambda_n(D) P_n \text{ for some } \Lambda_n(D)\in \operatorname{Mat}_N(\mathbb{C}), \text{ for all } n\in\mathbb{N}_0\Big\},
\end{equation} 
where $\{P_n\}_{n\in \mathbb{N}_0}$ is any sequence of matrix-valued orthogonal polynomials with respect to $W$ and as usual $\mathbb N_0=\{0,1,\dots  \}$.
We observe that the definition of $\mathcal D(W)$ depends only on the weight matrix $W$ and not on the particular sequence of orthogonal polynomials.

If  $\{P_n\}_{n\geq 0}$ is  the sequence of monic orthogonal polynomials associated to the weight $W$, and  $D=\sum_{i=0}^m \partial ^i F_i(x) \in \mathcal D(W)$, with $F_i(x)=\sum_{j=0}^i x^j F_j^i$, then  the eigenvalues $\Lambda_n$ are given by 
\begin{equation}\label{eigenvaluemonicos}
   \Lambda_n=\Lambda_n(D)=\sum_{i=0}^m [n]_i F_i^i, \qquad \text{for all } n\geq 0,
 \end{equation}
    where $[n]_i=n(n-1)\cdots (n-i+1)$, $[n]_0=1$.

The  formal adjoint on $\operatorname{Mat}_{N}(\Omega([[x]])$, denoted  by $\mbox{}^*$, is the unique involution extending Hermitian conjugation on $\operatorname{Mat}_{N}(\mathbb C[x])$ and mapping $\partial I$ to $-\partial I$. 
The {\em formal $W$-adjoint} of $ \mathfrak{D}\in \operatorname{Mat}_{N}(\Omega[x])$ is the differential operator $\mathfrak{D}^{\dagger} \in \operatorname{Mat}_{N}(\Omega[[x]])$ defined by
$$\mathfrak{D}^{\dagger}:= W(x)\mathfrak{D}^{\ast}W(x)^{-1},$$
where $\mathfrak{D}^{\ast}$ is the formal adjoint of $\mathfrak D$. 
An operator $\mathfrak{D}\in \operatorname{Mat}_{N}(\Omega[x])$ is called {\em $W$-adjointable} if there exists $\widetilde {\mathfrak{D}} \in \operatorname{Mat}_{N}(\Omega[x])$, such that
$$\langle P\cdot \mathfrak{D},Q\rangle=\langle P,Q\cdot \tilde{\mathfrak{D}}\rangle,$$
for all $P,Q\in \operatorname{Mat}_N(\mathbb{C}[x])$. Then we say that the operator $\tilde{\mathfrak D}$ is the $W$-adjoint of $\mathfrak D $.

\smallskip

In this paper, we work with weight matrices $W(x)$ that are 'good enough'. Specifically, these weights are such that for each integer $n \geq 0$, the $n$-th derivative $W^{(n)}(x)$ decreases exponentially at infinity, and there exists a scalar polynomial $p_n(x)$ such that $W^{(n)}(x)p_n(x)$ has finite moments. See \cite{CY18}*{Section 2.2}. With this in mind, the following proposition holds.

\begin{prop}[\cite{CY18}, Prop. 2.23] \label{adjuntas}  
If \( \mathfrak{D} \in \operatorname{Mat}_{N}(\Omega[x]) \) is \( W \)-adjointable, then its \( W \)-adjoint $\tilde{\mathfrak{D}}$ coincides with its formal $W$-adjoint; that is, \( \widetilde{\mathfrak{D}} = \mathfrak{D}^{\dagger} \).
\end{prop}

For a differential operator \( \mathfrak{D} = \sum_{j=0}^n \partial^j F_j \in \operatorname{Mat}_{N}(\Omega[x]) \), its formal \( W \)-adjoint \( \mathfrak{D}^\dagger \) is explicitly given by
\[
\mathfrak{D}^\dagger = \sum_{k=0}^n \partial^k G_k, \quad \text{where } G_k = \sum_{j = 0}^{n-k}(-1)^{n-j} \binom{n-j}{k} (W F_{n-j}^*)^{(n-k-j)} W^{-1}, \quad \text{for } 0 \leq k \leq n.
\]

\smallskip

The following result holds for the algebra $\mathcal{D}(W)$.

\begin{prop} \label{adjunta D(W)} 
    If $D \in \mathcal{D}(W)$, then  $D$ is $W$-adjointable. Moreover, the $W$-adjoint $D^{\dagger}$ also belongs to $\mathcal{D}(W)$.
\end{prop}
\begin{proof}
     The proof follows by combining Proposition \ref{adjuntas} with \cite{GT07}*{Cor. 4.5}.
\end{proof}

We now define $W$-symmetric operators, which play a fundamental role in the algebra \( \mathcal{D}(W) \).

\begin{definition}
A differential operator \( D \in \operatorname{Mat}_{N}(\Omega[x]) \) is said to be \( W \)-\emph{symmetric} if
\[
\langle P \cdot D, Q \rangle = \langle P, Q \cdot D \rangle, \quad \text{for all } P, Q \in \operatorname{Mat}_{N}(\mathbb{C}[x]).
\]
Equivalently, \( D \) is \( W \)-adjointable with \( D = D^\dagger \).
\end{definition}

\smallskip

The set \( \mathcal{S}(W) \) of all \( W \)-symmetric operators in \( \mathcal{D}(W) \) has the following key property
\[
\mathcal{D}(W) = \mathcal{S}(W) \oplus i \mathcal{S}(W),
\]
where \( i \) is the imaginary unit. This decomposition emphasizes the central role of \( W \)-symmetric operators in understanding the algebra \( \mathcal{D}(W) \).

\smallskip

To provide a more concrete setting, we observe that \( \mathcal{D}(W) \) is a subalgebra of the algebra
\begin{equation}\label{gradoscorrectos}
\mathbf{D} = \left\{ \sum_{j=0}^{m} \partial^{j} F_{j} \in \operatorname{Mat}_{N}(\Omega[x]) : \deg(F_j) \leq j \right\}.
\end{equation}

An important result, first established in \cite{GT07}, states that any \( W \)-symmetric operator in \( \mathbf{D} \) must belong to the algebra \( \mathcal{D}(W) \). 
We also have the following result, which connects the coefficients of a differential operator with its action on the sequence of monic orthogonal polynomials.

\begin{prop} \label{graleigenv}
Let \( D = \sum_{i=0}^m \partial^i F_i(x) \in \mathbf{D} \), where \( F_i(x) = \sum_{j=0}^i x^j F_j^i \) for some \( F_j^i \in \operatorname{Mat}_N(\mathbb{C}) \). Then, for any \( \nu \in \mathbb{C} \), the expression
\[
\Lambda_\nu(D) = \sum_{i=0}^m [\nu]_i F_i^i, \quad \text{where } [\nu]_i = \nu(\nu-1) \cdots (\nu - i + 1), \quad [\nu]_0 = 1,
\]
defines a multiplicative map on \( \mathbf{D} \) when restricted to square matrices, giving a representation of \( \mathbf{D} \) into \( \operatorname{Mat}_N(\mathbb{C}) \).
\end{prop}

\smallskip

\begin{remark}\label{Lambdanu}
This result extends naturally to matrix differential operators of general size:
Let \( D \in \operatorname{Mat}_{N \times M}(\Omega[x]) \) be an \( N \times M \)-matrix differential operator of the form
\[
D = \sum_{j=0}^m \partial^j F_j, \quad \text{with } \deg(F_j) \leq j.
\]
Then, for any \( \nu \in \mathbb{C} \), we define
\[
\Lambda_\nu(D) = \sum_{i=0}^m [\nu]_i F_i^i \in \operatorname{Mat}_{N \times M}(\mathbb{C}).
\]
In this setting, while \( \Lambda_\nu \) cannot be interpreted as a representation of an algebra (due to the differing sizes of the matrices), it still satisfies a {\em multiplicativity property} when the dimensions are compatible. Namely, for any differential operators \( E \in \operatorname{Mat}_{R \times N}(\Omega[x]) \) and \( D \in \operatorname{Mat}_{N \times M}(\Omega[x]) \), we have
$$
\Lambda_\nu(ED) = \Lambda_\nu(E)\Lambda_\nu(D). $$

This result will be particularly useful when studying the modules \( \mathcal{D}(W, \widetilde{W}) \) in later sections, where we deal with operators acting between spaces of orthogonal polynomials associated with weights of different sizes.
\end{remark}

\section{Darboux transformations}\label{Darboux-sect}
In this section, we recall the notion of Darboux transformations for weight matrices and their associated orthogonal polynomials. These transformations provide a method to relate different weight matrices through specific factorizations of differential operators. The results presented here establish the framework that will be used later to study new examples of matrix-valued orthogonal polynomials.

Throughout this section, we consider weight matrices $W$ and $\tilde{W}$ supported on the same interval. We denote $P_n=P_{n}(x)$ and $\widetilde P_n=\tilde{P}_{n}(x)$ their associated sequence of monic orthogonal polynomials, respectively.

We say that  $\mathcal{V} \in \operatorname{Mat}_{N}(\Omega [x])$, is a  \textbf{degree-preserving} operator if the degree of $P(x) \cdot \mathcal{V}$ is equal to the degree of $P(x)$ for all $P(x) \in \operatorname{Mat}_{N}(\mathbb{C})[x]$.

\begin{definition}\label{darb def}
    We say that $\tilde{W}$ is a Darboux transformation of $W$ if there exists a differential operator $D \in \mathcal{D}(W)$ that can be factorized as $D = \mathcal{V}\mathcal{N}$, with degree-preserving operators $\mathcal{V}, \, \mathcal{N} \in \operatorname{Mat}_{N}(\Omega[x])$, such that
    \begin{equation*}
            P_{n}(x) \cdot \mathcal{V}  = A_{n}\tilde{P}_{n}(x) , \qquad  \text{ for all } n \in \mathbb N_0, 
    \end{equation*}
for some matrices $A_{n} \in \operatorname{Mat}_{N}(\mathbb{C})$.
\end{definition}

\begin{remark}
   Historically, the term Darboux transformation refers broadly to the factorization of a differential operator without further restrictions. Here, we require the factors to be degree-preserving to ensure that  \(\widetilde W\) still yields a bispectral sequence of
orthogonal matrix polynomials.
\end{remark}

The definition of a Darboux transformation relies on the factorization of a differential operator in the algebra $\mathcal D(W)$ and the degree-preserving property. The following result provides a more practical characterization by reducing the problem to the existence of a single operator \( \mathcal{V} \) that transforms the orthogonal polynomials \( P_n \) of \( W \) into those of \( \tilde{W} \), with suitable invertibility conditions on the transformation matrices \( A_n \).

\begin{prop}\label{darboux-thm} (cf. \cite{BPZ24}*{Theorem  3.3})  A weight 
$\tilde {W}$ is a {Darboux transformation} of $W$ if and only if 
there exists a differential operator $\mathcal{V}\in \operatorname{Mat}_{N}(\Omega[x])$ such that 
    \begin{equation}\label{darboux-1}
            P_{n}(x) \cdot \mathcal{V}  = A_{n}\tilde {P}_{n}(x), \quad \text{ for all } n\in \mathbb{N}_{0}, 
    \end{equation}
where $A_{n} \in \operatorname{Mat}_{N}(\mathbb{C})$ invertible for all  
$n$.
\end{prop}

\begin{proof}
Let us observe that if $\tilde W $ is a Darboux transformation of $W$, then the matrices $A_n$ are invertible because $\mathcal V$ is a degree-preserving operator. 

Conversely, let us  consider \( \hat{P}_n := P_n \oplus \tilde{P}_n \) be the monic orthogonal polynomials corresponding to the weight \( \hat{W} := W \oplus \tilde{W} \). It is easy to see that 
the operator  $E  = \begin{psmallmatrix} 0 & \mathcal{V} \\ 0 & 0 \end{psmallmatrix} \in \mathcal D(\hat W)$ and $\Lambda_n(E) = \begin{psmallmatrix} 0 & A_n \\ 0 & 0 \end{psmallmatrix}.$ The adjoint of $E$ is 
is given by $E^\dagger = \begin{psmallmatrix} 0 & 0 \\ \mathcal{N} & 0 \end{psmallmatrix}$,  where $\mathcal{N} = \tilde{W} \mathcal{V}^* W^{-1}$ and 
\[
\Lambda_n(E^\dagger) = \begin{psmallmatrix} 0 & 0 \\ \tilde{A}_n & 0 \end{psmallmatrix}, \quad \text{with} \quad \tilde{A}_n = \|\tilde{P}_n\|_{\tilde{W}}^2 A_n^* \|P_n\|_W^{-2}.
\]
Thus, we have that 
$\tilde{P}_n \cdot \mathcal{N} = \tilde{A}_n P_n$, where  $\tilde{A}_n $ is invertible for all \( n \).
This proves that \( \mathcal{N} \) is also a degree-preserving operator.

Finally, since 
$P_n \cdot \mathcal{V} \mathcal{N} = A_n \tilde{A}_n P_n$,
we obtain 
\( D =\mathcal{V} \mathcal{N} \in \mathcal{D}(W) \). Therefore, \( \tilde{W} \) is a Darboux transformation of \( W \), as required.

\end{proof}

\begin{cor} \label{cons}
Let \( \tilde{W} \) be a Darboux transformation of \( W \), and let \( \mathcal{V} \) be a differential operator satisfying \eqref{darboux-1}. Define \( \mathcal{N} = \tilde{W} \mathcal{V}^* W^{-1} \). Then the following hold:
\begin{enumerate}[ i)]
    \item \( W \) is a Darboux transformation of \( \tilde{W} \).
\item The differential operators \( D = \mathcal{V} \mathcal{N} \in \mathcal{D}(W) \) and \( \tilde{D} = \mathcal{N} \mathcal{V} \in \mathcal{D}(\tilde{W}) \). Moreover, their eigenvalues \( \Lambda_n(D) \) and \( \Lambda_n(\tilde{D}) \) are nonsingular matrices for all \( n \geq 0 \).

    \item The operators
    \[
    \begin{pmatrix} 0 & \mathcal{V} \\ 0 & 0 \end{pmatrix} \quad \text{and} \quad \begin{pmatrix} 0 & 0 \\ \mathcal{N} & 0 \end{pmatrix}
    \]
    belong to the algebra \( \mathcal{D}(W \oplus \tilde{W}) \).
\end{enumerate}
\end{cor}

Darboux transformation defines an equivalence relation on the set of all weights of size $N$ supported in the same interval. This relation preserves direct sums of weights:
If $\tilde{W}_{1}$ and $\tilde{W}_{2}$ are Darboux transformations of $W_{1}$ and $W_{2}$, respectively, then $\tilde{W}_{1} \oplus \tilde{W}_{2}$ is a Darboux transformation of $W_{1} \oplus W_{2}$.

The Darboux transformation not only preserves the structure of direct sums of weights but also establishes a strong connection between the differential operator algebras associated with the original weight  $W$  and its transformed counterpart  $\tilde{W}$. In particular, the following result shows how these algebras interact through the operators  $\mathcal{V}$  and  $\mathcal{N}$  from the Darboux factorization.

\begin{prop}\label{alg debil}
Let \( W \) and \( \tilde{W} \) be weight matrices of size \( N \), such that \( \tilde{W} \) is a Darboux transformation of \( W \) with \( D = \mathcal{V} \mathcal{N} \in \mathcal{D}(W) \), as in Definition \ref{darb def}. Then, the  operator algebras satisfy the following inclusions
\[
\mathcal{V} \mathcal{D}(\tilde{W}) \mathcal{N} \subseteq \mathcal{D}(W)
\quad \text{and} \quad
\mathcal{N} \mathcal{D}(W) \mathcal{V} \subseteq \mathcal{D}(\tilde{W}).
\]
\end{prop}

\begin{proof}
We have that \( P_{n}(x) \cdot \mathcal{V} = A_{n} \tilde{P}_{n}(x) \) and \( \tilde{P}_{n}(x) \cdot \mathcal{N} = \tilde{A}_{n} P_{n}(x) \) for some invertible matrices \( A_{n} \) and \( \tilde{A}_{n} \).

Let \( \mathcal{A} \in \mathcal{D}(\tilde{W}) \) and fix \( n \geq 0 \). To prove that \( \mathcal{V} \mathcal{A} \mathcal{N} \in \mathcal{D}(W) \), it is enough to show that \( P_{n}(x) \) is an eigenfunction of \( \mathcal{V} \mathcal{A} \mathcal{N} \). We compute
\[
P_{n}(x) \cdot \mathcal{V} \mathcal{A} \mathcal{N} 
= \big( P_{n}(x) \cdot \mathcal{V} \big) \cdot \mathcal{A} \mathcal{N} 
= \big( A_{n} \tilde{P}_{n}(x) \big) \cdot \mathcal{A} \mathcal{N}.
\]
Since \( \mathcal{A} \in \mathcal{D}(\tilde{W}) \), we know that \( \tilde{P}_{n}(x) \cdot \mathcal{A} = \Lambda_{n}(\mathcal{A}) \tilde{P}_{n}(x) \) for some matrix \( \Lambda_{n}(\mathcal{A}) \). Substituting this, we obtain
\[
P_{n}(x) \cdot \mathcal{V} \mathcal{A} \mathcal{N} 
= A_{n} \Lambda_{n}(\mathcal{A}) \big( \tilde{P}_{n}(x) \cdot \mathcal{N} \big)
= A_{n} \Lambda_{n}(\mathcal{A}) \tilde{A}_{n} P_{n}(x).
\]
Thus, \( P_{n}(x) \) is indeed an eigenfunction of \( \mathcal{V} \mathcal{A} \mathcal{N} \), and we conclude that \( \mathcal{V} \mathcal{A} \mathcal{N} \in \mathcal{D}(W) \).

A similar argument shows that \( \mathcal{N} \mathcal{D}(W) \mathcal{V} \subseteq \mathcal{D}(\tilde{W}) \), completing the proof.
\end{proof}

One of the main challenges in identifying when two weights are related by a Darboux transformation lies in determining their corresponding sequences of orthogonal polynomials. In practice, these sequences are often difficult to compute explicitly. The following proposition, originally proposed and proven in \cite{BPstrong}, provides an alternative approach: it gives sufficient conditions to establish that a weight $\tilde{W}$ is a Darboux transformation of $W$ without requiring prior knowledge of their orthogonal polynomials. This result is particularly useful in cases where only partial information about the differential operators or the weights is available.

\begin{prop}\label{strong1}
     If there exists a degree-preserving differential operator \( \mathcal{V} \) such that \( D = \mathcal{V}\mathcal{N} \in \mathcal{D}(W) \), with \( \mathcal{N} = \tilde{W}\mathcal{V}^{\ast}W^{-1} \), and
    \[
    \langle P\cdot \mathcal{V}, Q \rangle_{\tilde{W}} = \langle P, Q\cdot \mathcal{N} \rangle_{W} \qquad \text{ for all } P,Q \in \operatorname{Mat}_{N}(\mathbb{C}[x]),
    \]
    then
    \begin{enumerate} [ (i)]
        \item \( \tilde{W} \) is a Darboux transformation of \( W \).

        \item The symmetric operators in the algebra \( \mathcal{D}(\widetilde{W}) \) satisfy
        \[
        \mathcal{V}\mathcal{S}(\tilde{W})\mathcal{N} \subseteq \mathcal{S}(W).
        \]

        \item If \( \widetilde{D} \in \mathcal{D}(\widetilde{W}) \) is of even (odd) order, and the leading coefficient of \( \mathcal{V} \) is nonsingular, then \( D = \mathcal{V}\widetilde{D}\mathcal{N} \in \mathcal{D}(W) \) is of even (odd) order.
    \end{enumerate}
\end{prop}
\begin{proof}
{\it (i)}  Let \( P_{n} \) and \( \tilde{P}_{n} \) be the sequences of monic orthogonal polynomials associated with \( W \) and \( \tilde{W} \), respectively. Since \( \mathcal{V} \) is degree-preserving, \( P_{n}(x) \cdot \mathcal{V} \) is a polynomial of degree \( n \) with nonsingular leading coefficient. For \( n \neq m \), we have
\[
\langle P_{n} \cdot \mathcal{V}, P_{m} \cdot \mathcal{V} \rangle_{\tilde{W}} = \langle P_{n}, P_{m} \cdot \mathcal{V} \mathcal{N} \rangle_{W} = \langle P_{n}, P_{m} \rangle_{W} \Lambda_{n}(D)^{\ast} = 0,
\]
where \( D = \mathcal{V} \mathcal{N} \in \mathcal{D}(W) \) and \( \Lambda_n(D) \) is nonsingular. Therefore, \( P_{n}(x) \cdot \mathcal{V} \) is a sequence of orthogonal polynomials with respect to \( \tilde{W} \), and hence there exist nonsingular matrices \( A_n \) such that
\[
P_{n}(x) \cdot \mathcal{V} = A_{n} \tilde{P}_{n}(x).
\]
Thus, \( \tilde{W} \) is a Darboux transformation of \( W \).

{\it (ii)}  Let \( \mathcal{A} \in \mathcal{S}(\tilde{W}) \), i.e., \( \mathcal{A} = \mathcal{A}^\dagger = \tilde{W} \mathcal{A}^* \tilde{W}^{-1} \). By Proposition \ref{alg debil}, we have that \( \mathcal{B} = \mathcal{V} \mathcal{A} \mathcal{N} \in \mathcal{D}(W) \). To show that \( \mathcal{B} \) is \( W \)-symmetric, we compute its adjoint
\[
\mathcal{B}^\dagger = W \mathcal{B}^* W^{-1} = W \mathcal{N}^* \mathcal{A}^* \mathcal{V}^* W^{-1}.
\]
Using the definitions of \( \mathcal{A} \) and \( \mathcal{N} \), we find that 
$\mathcal{B}^\dagger = \mathcal{V} \mathcal{A} \mathcal{N} = \mathcal{B}$.
Thus, \( \mathcal{B} \in \mathcal{S}(W) \).

{\it (iii)}  Let \( \widetilde{D} \in \mathcal{D}(\tilde{W}) \), and define \( \mathcal{D} = \mathcal{V} \widetilde{D} \mathcal{N} \). Since \( \mathcal{V} \) has nonsingular leading coefficient, it follows that \( \mathcal{N} = \tilde{W} \mathcal{V}^\ast W^{-1} \) also has nonsingular leading coefficient, and \( \operatorname{Ord}(\mathcal{V}) = \operatorname{Ord}(\mathcal{N}) = \ell \). Hence
\[
\operatorname{Ord}(\mathcal{D}) = \operatorname{Ord}(\mathcal{V}) + \operatorname{Ord}(\widetilde{D}) + \operatorname{Ord}(\mathcal{N}) = 2\ell + \operatorname{Ord}(\widetilde{D}).
\]
Therefore
\( \operatorname{Ord}(\mathcal{D}) \) has the same parity as \( \operatorname{Ord}(\widetilde{D}) \). 
\end{proof}

\begin{remark}
Proposition \ref{strong1} provides a practical tool to obtain a sequence of orthogonal polynomials for \( \tilde{W} \) when the sequence for \( W \) is already known. In particular, if \( \{P_n\} \) denotes the sequence of monic orthogonal polynomials with respect to \( W \), then the sequence \( \{P_n \cdot \mathcal{V}\} \), where \( \mathcal{V} \) is the degree-preserving operator in the factorization \( D = \mathcal{V}\mathcal{N} \), forms a sequence of orthogonal polynomials with respect to \( \tilde{W} \).
Moreover,
    $$P_{n}(x) \cdot \mathcal{V} = \Lambda_{n}(\mathcal{V})\tilde{P}_{n}(x),$$
    where  $\Lambda_n(\mathcal{V})$  are the eigenvalues associated with  $\mathcal{V}$, defined in Proposition \ref{graleigenv}.
\end{remark}

To illustrate the concept of Darboux transformations in the scalar case, we consider the classical scalar weights on the real line. These weights are characterized by the existence of a second-order differential operator $\delta$ in their algebra \( \mathcal{D}(w) \), making them solutions to the classical Bochner problem in the scalar case.
Up to an affine change of variables, the only such weights are the Hermite, Laguerre, and Jacobi weights, whose corresponding differential operators are given in Table \ref{weights-table}. These operators play a crucial role in the theory of orthogonal polynomials and will serve as key components in the analysis of the matrix-valued case.

\begin{table}[H] 
\begin{center} \small
    \begin{tabular}{|c|c|c|c|}
     \hline    Type & Differential operator  & Weight   & Support    \\    \hline &&& \\ 
     (shifted) Hermite   &   $\delta =  \partial^{2} -2 \partial (x-b)$ & $w_b(x)=e^{-x^2+2bx}$ & $(-\infty, \infty)$   \\ 
     \hline & & & \\   
     Laguerre    &  $\delta =  \partial^{2}x  + \partial(\alpha+1-x)$   & $w_\alpha(x)= e^{-x} x^\alpha$ & $ (0,\infty )$ \\ 
     \hline  & & & \\
     Jacobi  & $\delta = \partial^{2} (1-x^2)  + \partial (\beta-\alpha-x(\alpha+\beta+2))$ & $w_{\alpha,\beta}(x) = (1-x)^\alpha (1+x)^\beta $ & $(-1,1)$ \\ & & & \\\hline
    \end{tabular}
    \caption{Classical scalar weights. } 
    \label{weights-table}
\end{center} 
\end{table}

\begin{remark} 
    In addition to the classical Hermite weights, we include the shifted Hermite weights \( w_b(x) \) in this analysis. These weights arise naturally in the construction of several matrix-valued weights that solve the Bochner problem in the matrix case. 
\end{remark}

The classification of Darboux-equivalent classes for these weights was established in \cite{BP23b} and it will play a key role in understanding matrix weights that arise as direct sums or Darboux transformations of classical scalar weights.

\begin{prop}\label{scalar-darboux}
The Darboux-equivalent classes of the classical scalar weights are as follows:
\begin{itemize}
    \item A shifted Hermite weight \( w_b(x) \) is not Darboux-equivalent to \( w_a(x) \) if \( a \neq b \).
    \item A Laguerre weight \( w_\alpha(x) \) is Darboux-equivalent to \( w_{\widetilde{\alpha}}(x) \) if and only if \( \alpha - \widetilde{\alpha} \in \mathbb{Z} \).
    \item A Jacobi weight \( w_{\alpha,\beta}(x) \) is Darboux-equivalent to \( w_{\widetilde{\alpha},\widetilde{\beta}}(x) \) if and only if \( \alpha + \beta = \widetilde{\alpha} + \widetilde{\beta} \) and \( \alpha - \widetilde{\alpha} \in \mathbb{Z} \).
\end{itemize}
\end{prop}

\section{Centers of algebras under Darboux transformations}\label{sec-centers}

In this section, we focus on the center \( \mathcal{Z}(W) \) of the algebra \( \mathcal{D}(W) \) associated with a matrix weight \( W \) and explore its behavior under Darboux transformations. While the previous sections provided a survey of known results on Darboux transformations and their connection to differential operator algebras, here we develop new results that highlight the interplay between the center \( \mathcal{Z}(W) \) and its counterpart for a Darboux-equivalent weight \( \tilde{W} \).

We begin by analyzing the case where the operator \( D = \mathcal{V} \mathcal{N} \), as in Definition \ref{darb def}, belongs to the center \( \mathcal{Z}(W) \). This condition establishes a direct relationship between \( \mathcal{Z}(W) \) and \( \mathcal{Z}(\tilde{W}) \).

\begin{prop}\label{cent}
    Let \( \tilde{W} \) be a Darboux transformation of \( W \), with $\mathcal{V}$ and $\mathcal{N}$ as in Definition \ref{darb def}. If $D=\mathcal{V}\mathcal{N}$ belongs to the center of $\mathcal{D}(W)$ then
    \[
        \mathcal{N} \mathcal{Z}(W) \mathcal{V} \subseteq \mathcal{Z}(\tilde{W}) \quad \text{and} \quad \mathcal{V} \mathcal{Z}(\tilde{W}) \mathcal{N} \subseteq \mathcal{Z}(W).
        \]
        In particular, $\tilde{D} = \mathcal{N} \mathcal{V}$ belongs to $\mathcal{Z}(\tilde{W}) $.
\end{prop}
\begin{proof}
    Let \( B \in \mathcal{D}(\tilde{W}) \). By Proposition \ref{alg debil}, the operator \( \mathcal{V} B \mathcal{N} \in \mathcal{D}(W) \) and commutes with \( \mathcal{V} \mathcal{N} \in \mathcal{Z}(W) \). Therefore, for any \( Z \in \mathcal{Z}(W) \), we have
    \[
    \mathcal{V} \mathcal{N} Z \mathcal{V} B \mathcal{N} = \mathcal{V} B \mathcal{N} Z \mathcal{V} \mathcal{N}.
    \]
    Since \( \mathcal{V} \) and \( \mathcal{N} \) are not zero divisors, we can cancel \( \mathcal{V} \) and \( \mathcal{N} \) on both sides, which implies
    \[
    (\mathcal{N} Z \mathcal{V}) B = B (\mathcal{N} Z \mathcal{V}).
    \]
    This shows that \( \mathcal{N} Z \mathcal{V} \) commutes with all operators \( B \in \mathcal{D}(\tilde{W}) \), and therefore \( \mathcal{N} Z \mathcal{V} \in \mathcal{Z}(\tilde{W}) \). In particular, by taking \( Z = I \), we obtain \( \mathcal{N} \mathcal{V} \in \mathcal{Z}(\tilde{W}) \).
\end{proof}

In the following proposition, we establish one of our main results. We show that if a weight \( W \) is a Darboux transformation of a direct sum of classical scalar weights \( \tilde{W} \), then the Darboux transformation can always be obtained by factoring an operator \( D \) in the center \( \mathcal{Z}(\tilde{W}) \). This result, combined with the previous proposition, demonstrates an intrinsic connection between the centers \( \mathcal{Z}(W) \) and \( \mathcal{Z}(\tilde{W}) \). Leveraging the explicit knowledge of \( \mathcal{Z}(\tilde{W}) \), which we will develop in Section \ref{sect-escalar}, we can determine the center \( \mathcal{Z}(W) \) in the current setting.

\begin{thm}\label{darb center}
    Let \( W \) be a weight matrix that is Darboux-equivalent to \( \tilde{W} =  w_1 \oplus \cdots \oplus w_N\), where \( w_i \) are classical scalar weights. Then \( W \) can be obtained as a Darboux transformation of \( \tilde{W} \) by factoring a differential operator in the center of its algebra \( \mathcal{D}(\tilde{W}) \).
\end{thm}

\begin{proof}
Since the Darboux transformation defines an equivalence relation, we can assume without loss of generality that \( W \) is a Darboux transformation of 
\[
\tilde{W} = w_{1} I_{r_1} \oplus w_{2} I_{r_2} \oplus \cdots \oplus w_{m} I_{r_m},
\]
where \( w_{i} \not\sim w_{j} \) for \( i \neq j \), and \( r_1 + r_2 + \dots + r_m = N \).

Let \( D \in \mathcal{D}(W) \) with $D = \mathcal{V}\mathcal{N}$ as in the Definition \ref{darb def}. From Theorem \ref{eigenv-directsum} the eigenvalue of $D$ are of the form
\[
\Lambda_n(D) = \begin{psmallmatrix} \Lambda_n(D_1) & & \\ & \ddots & \\ & & \Lambda_n(D_m) \end{psmallmatrix},
\]
where \( \Lambda_n(D_i) \) belongs to \( \operatorname{Mat}_{r_i}(\mathbb{C}[\Lambda_n(\delta_i)]) \), with \( \delta_i \) being the second-order differential operator associated with \( w_i \).

For each \( i \), we note that \( \Lambda_n(D_i) \) is invertible, and therefore there exists a polynomial \( p_i \in \mathbb{C}[x] \) such that \( p_i(\Lambda_n(\delta_i)) \Lambda_n(D_i)^{-1} \) belongs to \( \operatorname{Mat}_{r_i}(\mathbb{C}[\Lambda_n(\delta_i)]) \). By Theorem \ref{eigenv-directsum}  we have that there exists a differential operator \( \mathcal{A} \in \mathcal{D}(\tilde{W}) \) such that 
\[
\Lambda_n(\mathcal{A}) \Lambda_n(D) = p_1(\Lambda_n(\delta_1)) I_{r_1} \oplus \cdots \oplus p_m(\Lambda_n(\delta_m)) I_{r_m}.
\]
Thus, by Theorem \ref{cent direct} the operator \( \mathcal{A}D \) belongs to the center \( \mathcal{Z}(\tilde{W}) \), and the Darboux transformation of \( \tilde{W} \) can be obtained by factoring \( \mathcal{A}D \). We have $\mathcal{A}D = (\mathcal{A}\mathcal{V})\mathcal{N}$, with
\[
P_n(x) \cdot \mathcal{A} \mathcal{V} = \left( p_1(\Lambda_n(\delta_1)) I_{r_1} \oplus \cdots \oplus p_m(\Lambda_n(\delta_m)) I_{r_m} \right) A_n \tilde{P}_n(x),
\]
where \( \mathcal{V} \) is the degree-preserving operator defining the Darboux transformation. This completes the proof.
\end{proof}

To conclude this section, we obtain an important structural property of the algebras  $\mathcal{D}(W)$  and  $\mathcal{D}(\tilde{W})$  when  $W$  and  $\tilde{W}$  are Darboux-equivalent. In particular, we show that the commutativity of the algebra  $\mathcal{D}(W)$  is preserved under Darboux transformations. This result further emphasizes the strong connection between the differential operator algebras associated with Darboux-equivalent weights.

\begin{thm}\label{conmutativa}
If $\tilde{W}$ is a Darboux transformation of $W$ then $\mathcal{D}(W)$ is a commutative algebra if and only if $\mathcal{D}(\tilde{W})$ is a commutative algebra.
\end{thm}
\begin{proof}
    Let \( \mathcal{A}_{1}, \mathcal{A}_{2} \in \mathcal{D}(\tilde{W}) \). Since \( \tilde{W} \) is a Darboux transformation of \( W \), we know from Proposition \ref{cent} that \( \mathcal{V} \mathcal{A}_{1} \mathcal{N} \in \mathcal{D}(W) \) and \( \mathcal{V} \mathcal{A}_{2} \mathcal{N} \in \mathcal{D}(W) \). 
    
 Assume \( \mathcal{D}(\tilde{W}) \) is commutative. Then, we have
    \[
    (\mathcal{V} \mathcal{A}_{1} \mathcal{N})(\mathcal{V} \mathcal{A}_{2} \mathcal{N}) 
    = \mathcal{V} (\mathcal{A}_{1} \mathcal{A}_{2}) \mathcal{N}
    = \mathcal{V} (\mathcal{A}_{2} \mathcal{A}_{1}) \mathcal{N}
    = (\mathcal{V} \mathcal{A}_{2} \mathcal{N})(\mathcal{V} \mathcal{A}_{1} \mathcal{N}).
    \]
    Hence, the algebra \( \mathcal{D}(W) \) is commutative.

    Conversely, assume \( \mathcal{D}(W) \) is commutative. For any \( \mathcal{A}_{1}, \mathcal{A}_{2} \in \mathcal{D}(\tilde{W}) \), consider the operators \( \mathcal{V} \mathcal{A}_{1} \mathcal{N} \) and \( \mathcal{V} \mathcal{A}_{2} \mathcal{N} \), which belong to \( \mathcal{D}(W) \). Since \( \mathcal{D}(W) \) is commutative, we have
    \[
    (\mathcal{V} \mathcal{A}_{1} \mathcal{N})(\mathcal{V} \mathcal{A}_{2} \mathcal{N}) 
    = (\mathcal{V} \mathcal{A}_{2} \mathcal{N})(\mathcal{V} \mathcal{A}_{1} \mathcal{N}).
    \]
    Using the fact that \( \mathcal{V} \) and \( \mathcal{N} \) are not zero divisors, it follows that \( \mathcal{A}_{1} \mathcal{A}_{2} = \mathcal{A}_{2} \mathcal{A}_{1} \), proving that \( \mathcal{D}(\tilde{W}) \) is commutative.
\end{proof}

\section{Eigenvalue algebras and modules of differential operators}\label{sect-eigenvalues}

In this section, we analyze two central structures associated with matrix weights: the eigenvalues algebra $\Lambda(W)$ and the modules $\mathcal D(W, \tilde{W})$. 
The eigenvalues algebra $\Lambda(W)$ captures the structure of $D(W)$ by encoding the action of its differential operators on the sequence of monic orthogonal polynomials with respect to $W$. 
 On the other hand, the modules $\mathcal D(W, \widetilde{W})$ describe differential operators that map one sequence of orthogonal polynomials, associated with a weight $W$, into another sequence associated with a potentially different weight $\tilde{W}$. Together, these structures offer complementary perspectives on the interplay between differential operators, allowing us to explore both the intrinsic algebraic properties of $D(W)$ and the relationships between distinct weights.

Associated with any $D \in \mathcal{D}(W)$, we define the function 
\[
\Lambda(D) : n \longmapsto \Lambda_n(D),
\]
where $\Lambda_n(D)$ is the eigenvalue associated with the $n$-th monic orthogonal polynomial. Intrinsically, we defined it in terms of the coefficients of the differential operator by 
\[
\Lambda_n(D) = \Lambda\Big(\sum_{j=0}^m \partial^j F_j(x)\Big) 
= \sum_{j=0}^m [n]_j F_j^j,
\]
where $[n]_j = n(n-1)\cdots(n-j+1)$ for $j \geq 1$ and $[n]_0 = 1$.
Since $\Lambda_n(D)$ is a polynomial function of $n$, the map $\Lambda(D)$ belongs to $\operatorname{Mat}_N(\mathbb{C}[n])$, the space of $N \times N$ matrices with polynomial entries in the variable $n$. Alternatively, we can identify $\Lambda(D)$ with the sequence $$\left( \Lambda_0(D), \Lambda_1(D), \Lambda_2(D), \dots \right).$$

This sequence of eigenvalues $\Lambda_n(D)$ significantly simplifies many calculations related to the algebra $D(W)$. By focusing on the eigenvalue maps, we obtain a clearer and more efficient framework for analyzing the algebraic structure of differential operators associated with matrix weights.

\begin{prop}
The map 
\[
\Lambda : \mathcal{D}(W) \longrightarrow \operatorname{Mat}_N(\mathbb{C}[n]), \qquad D \mapsto \Lambda(D),
\]
is an injective homomorphism of algebras, meaning that 
\[
\Lambda(D_1 D_2) = \Lambda(D_1) \Lambda(D_2),
\]
for all $D_1, D_2 \in \mathcal{D}(W)$.
\end{prop}

\begin{remark}
The map $\Lambda$ is a faithful representation of the algebra $\mathcal{D}(W)$ in the infinite-dimensional Hilbert space $\operatorname{Mat}_N(\mathbb{C}[n])$.
\end{remark}

\begin{proof}
From Proposition \ref{graleigenv}, the mapping $D \mapsto \Lambda_n(D)$ is a representation of $\mathcal{D}(W)$ in $\operatorname{Mat}_N(\mathbb{C})$ for any $n \in \mathbb{N}_0$. Therefore, 
\[
\Lambda_n(D_1 D_2) = \Lambda_n(D_1) \Lambda_n(D_2),
\]
for all $D_1, D_2 \in \mathcal{D}(W)$.

Finally, suppose that $D = \sum_{j=0}^m \partial^j F_j(x)$ satisfies $\Lambda(D) = 0$. Then $\Lambda_n(D) = 0$ for all $n \in \mathbb{N}_0$. In particular, $F_0 = \Lambda_0(D) = 0$, and proceeding inductively, we find that $F_j = 0$ for all $j$. Hence, $D = 0$, and the map $\Lambda$ is injective.
\end{proof}

Using the homomorphism $\Lambda$, we can naturally identify $\mathcal{D}(W)$ as a subalgebra of $\operatorname{Mat}_N(\mathbb{C}[n])$, the algebra of $N \times N$ matrices with entries in $\mathbb{C}[n]$.

\begin{definition}
Given a weight matrix $W$ of size $N$, we define the eigenvalues algebra associated with $W$ as
\[
\Lambda(W) = \Lambda \left( \mathcal{D}(W) \right) = \big\{ \Lambda(D) : D \in \mathcal{D}(W) \big\} \subset \operatorname{Mat}_N(\mathbb{C}[n]).
\]
\end{definition}

The injectivity of \( \Lambda \) is a crucial property, as it ensures that the algebra \( \mathcal{D}(W) \) of differential operators is faithfully represented within \( \operatorname{Mat}_N(\mathbb{C}[n]) \). In particular, it establishes a one-to-one correspondence between elements of \( \mathcal{D}(W) \) and their images under \( \Lambda \).

While \( \mathcal{D}(W) \) and its eigenvalue algebra \( \Lambda(W) \) are isomorphic, and thus share the same algebraic structure and properties, it is often advantageous to work directly with \( \Lambda(W) \). The eigenvalue algebra offers a computationally simpler framework, as it consists of matrices with polynomial entries in \( \mathbb{C}[n] \). 
This alternative perspective allows us to efficiently derive strong results about the algebra \( \mathcal{D}(W) \) and its center \( \mathcal{Z}(W) \), particularly in the context of Darboux transformations. See, for instance, Theorem \ref{darb center}.

\begin{cor}
A differential operator $D \in \mathcal{D}(W)$ belongs to the center of $\mathcal{D}(W)$ if and only if $\Lambda(D)$ lies in the center of $\Lambda(W)$.
\end{cor}

For classical scalar weights such as Hermite, Laguerre, and Jacobi, the monic orthogonal polynomials $p_n$ are eigenfunctions of the differential operator $\delta$, given in Table \ref{weights-table}. Their eigenvalues are:
\begin{equation}\label{class-eigenv}
\Lambda_n(\delta) =
\begin{cases}
-2n, & \text{for Hermite weights}, \\
-n, & \text{for Laguerre weights}, \\
-n(n+\alpha+\beta+1), & \text{for Jacobi weights $w_{\alpha,\beta}$}.
\end{cases}
\end{equation}

Using these eigenvalues, we can explicitly describe the eigenvalues algebra $\Lambda(w)$ for classical weights. Since $\mathcal{D}(w)$ is a polynomial algebra generated by the operator $\delta$, the algebra $\Lambda(w)$ is:
\begin{equation}\label{AE-escalar}
\Lambda(w) = \mathbb{C}[\Lambda_n(\delta)] =
\begin{cases}
\mathbb{C}[n], & \text{if $w$ is a Laguerre or Hermite weight}, \\
\mathbb{C}[n(n+\alpha+\beta+1)], & \text{if $w = w_{\alpha, \beta}$ is a Jacobi weight}.
\end{cases}
\end{equation}

\

\smallskip
Let $W$ and $\widetilde{W}$ be weight matrices of sizes $N_1$ and $N_2$, respectively, with $P_n$ and $\widetilde{P}n$ denoting their sequences of monic orthogonal polynomials. We recover the modules $\mathcal{D}(W, \widetilde{W})$ from \cite{TZ18}, which extend the framework of $\mathcal{D}(W)$ to study connections between different weight matrices, even when they have different sizes. We define
 \begin{equation}\label{modulos Ww}
            \mathcal{D}(W,\tilde  W)  = \left\{ \mathcal T \in \operatorname{Mat}_{N_{1} \times N_{2}}(\Omega[x]) \, : \, P_{n}(x) \cdot \mathcal T = A_n \tilde  P_{n}(x), \text{ with } A_n \in \operatorname{Mat}_{N_{1}\times N_{2}}(\mathbb{C}) \right \}.
    \end{equation}

Observe that $\mathcal{D}(W, \widetilde{W})$ is naturally both a left $\mathcal{D}(W)$-module and a right $\mathcal{D}(\widetilde{W})$-module. These properties show the role of $\mathcal{D}(W, \widetilde{W})$ in linking the algebraic structures associated with $W$ and $\widetilde{W}$. 

The study of these modules has been developed extensively, with the general case discussed in \cite{BPZ24} and scalar cases treated in \cite{BP23b}. These results set the stage for the following straightforward proposition, which explores embeddings and adjoints of operators in $\mathcal{D}(W, \widetilde{W})$.

\begin{prop}\label{basicprop}
    Let $W$ and $\widetilde{W}$ be weight matrices of sizes $N_1$ and $N_2$, respectively. Let $\mathcal{T} \in \mathcal{D}(W, \widetilde{W})$ and $\mathcal{S} \in \mathcal{D}(\widetilde{W}, W)$. Then we have
    \begin{enumerate}
        \item[i)]  $\mathcal T \mathcal S\in \mathcal D(W)$ and  $\mathcal S \mathcal T\in \mathcal D(\widetilde W)$.
        
        \item[ ii)]The 
        modules $\mathcal D(W, \tilde W) $ and $\mathcal D( \tilde  W,  W) $  can be embedded into $\mathcal{D}(W \oplus \widetilde{W}):$ 
        \[
        T = \begin{pmatrix} 0 & \mathcal{T} \\ 0 & 0 \end{pmatrix} \in \mathcal{D}(W \oplus \widetilde{W}) \qquad  \text{ and }\qquad S = \begin{pmatrix} 0 & 0 \\ \mathcal{S} & 0 \end{pmatrix} \in \mathcal{D}(W \oplus \widetilde{W}).
        \]

        \item[iii)] The adjoint operator
        $\mathcal{T}^\dagger := \widetilde{W}(x) \mathcal{T}^* W^{-1}(x)$ belongs to $\mathcal{D}(\widetilde{W}, W)$.
    \end{enumerate}
\end{prop}

For any $\mathcal{T} \in \mathcal{D}(W, \widetilde{W})$, we define the map
\[
\Lambda(\mathcal{T}) : n \mapsto \Lambda_n(\mathcal{T}),
\]
where $\Lambda_n(\mathcal{T})$ is expressed in terms of the coefficients of $\mathcal{T}$ as
\[
\Lambda_n(\mathcal{T}) = \sum_{i=0}^m [n]_i F_i^i \in \operatorname{Mat}_{N_1 \times N_2}(\mathbb{C}),
\]
with $[n]_i = n(n-1)\cdots(n-i+1)$ for $i \geq 1$ and $[n]_0 = 1$. This map encodes how $\mathcal{T}$ interacts with the sequence of monic orthogonal polynomials associated with $W$ and $\widetilde{W}$:  For each $n$, we have
\[
P_{n}(x) \cdot \mathcal{T} = \Lambda_n(\mathcal{T}) \widetilde{P}_{n}(x),
\]
where $P_n(x)$ and $\widetilde{P}_n(x)$ denote the monic orthogonal polynomials of $W$ and $\widetilde{W}$, respectively.

The map $\Lambda$ is multiplicative, preserving the algebraic structure of the modules, in the sense that
\[
\Lambda(E_1 \mathcal{T} E_2) = \Lambda(E_1) \Lambda(\mathcal{T}) \Lambda(E_2),
\]
for all $E_1 \in \mathcal{D}(W)$ and $E_2 \in \mathcal{D}(\widetilde{W})$ (see Remark \ref{Lambdanu}). This property ensures compatibility between the action of $\mathcal{D}(W)$ and $\mathcal{D}(\widetilde{W})$ and the corresponding eigenvalue maps $\Lambda(E_1)$, $\Lambda(\mathcal{T})$, and $\Lambda(E_2)$, even when the matrices involved are not square.

\smallskip

Having introduced the general framework of modules $\mathcal{D}(W, \widetilde{W})$, we now focus on the classical scalar weights of Hermite, Laguerre, and Jacobi.
The structure of the modules $\mathcal{D}(w, \tilde{w})$ for these classical weights was thoroughly analyzed in \cite{BP23b}. In these cases, the modules are cyclic, meaning there exists a differential operator $\mathcal{T} = \mathcal{T}_{w, \tilde{w}}$ with polynomial coefficients that generates the module. Explicitly, we have
\[
\mathcal{D}(w, \tilde{w}) = \mathcal{T} \cdot \mathcal{D}(\tilde{w}) = \mathcal{D}(w) \cdot \mathcal{T}.
\]

The algebra $\mathcal{D}(w)$, in turn, is a polynomial algebra generated by the second-order differential operator $\delta$. The explicit expressions for the operators $\mathcal{T}$, which uniquely define the structure of these modules, are presented in the following theorem.

\begin{thm} \label{T-operators}
Let $w$ and $\tilde{w}$ be classical scalar weights. The differential operator $\mathcal{T} = \mathcal{T}_{w, \tilde{w}}$, such that
\[
\mathcal{D}(w, \tilde{w}) = \mathcal{T} \cdot \mathcal{D}(\tilde{w}),
\]
is given explicitly as follows:
\begin{enumerate}
    \item [(i)] $\mathcal{T}_{w, \tilde{w}} = 1$ if $w = \tilde{w}$, and $\mathcal{T}_{w, \tilde{w}} = 0$ if $w$ is not a Darboux transformation of $\tilde{w}$.
    \item [(ii)] For Laguerre weights $w_\alpha$, and $k \in \mathbb{N}$, we have
    \[
    \mathcal{T}_{w_\alpha, w_{\alpha+k}} = (-\partial + 1)^k, \qquad 
    \mathcal{T}_{w_{\alpha+k}, w_\alpha} = (\partial x + \alpha + k) \cdots (\partial x + \alpha + 1).
    \]
    \item [(iii)] For Jacobi weights $w_{\alpha,\beta}$, and $k \in \mathbb{N}$, we have
    \begin{align*}
    \mathcal{T}_{w_{\alpha,\beta+k}, w_{\alpha+k,\beta}} & = \big(\partial(1+x) + \beta + k\big) \cdots \big(\partial(1+x) + \beta + 1\big), \\
    \mathcal{T}_{w_{\alpha+k,\beta}, w_{\alpha,\beta+k}} & = \big(\partial(x-1) + \alpha + k\big) \cdots \big(\partial(x-1) + \alpha + 1\big).
    \end{align*}
\end{enumerate}
\end{thm}

\section{The algebra $\mathcal D(W)$ for a direct sum of classical weights}\label{sect-escalar}

In this section, we analyze the case where the weight $W$ is a direct sum of classical weights. Although such weights are not irreducible, the associated algebra of differential operators, $\mathcal{D}(W)$, exhibits a structure that goes beyond being a mere direct sum of the corresponding classical algebras. In particular, we identify operators in $\mathcal{D}(W)$ that cannot be expressed as direct sums of classical operators, highlighting the non-triviality of this scenario. This detailed exploration 
plays a crucial role in understanding more complex cases. The insights obtained here will later be applied to study irreducible weights that arise as Darboux transformations of classical weights.

\smallskip
 We consider weights that are direct sums of classical scalar weights of the same type (shifted Hermite, Laguerre, or Jacobi)
\[
W = w_{1} \oplus w_{2} \oplus \cdots \oplus w_{N}.
\]

The following result provides a full description of the algebra of differential operators in this setting.  
It is a direct consequence of  Theorems 4.5 and 4.6, and for this reason, we omit the proof.

\begin{thm}\label{algdirect}
Let $W = w_{1} \oplus w_{2} \oplus \cdots \oplus w_{N}$ be a direct sum of classical scalar weights of the same type (shifted Hermite, Laguerre, or Jacobi). Then, the algebra of differential operators associated with $W$ is given by
\[
\mathcal{D}(w_{1} \oplus \cdots \oplus w_{N}) = \sum_{i,j=1}^{N} \mathcal{T}_{w_i,w_j} \cdot \mathcal{D}(w_{j}) \, E_{i,j},
\]
where the differential operators $\mathcal{T}_{w_i,w_j}$ are as defined in Theorem \ref{T-operators}.
\end{thm}

As a direct application of Theorem \ref{algdirect}, we can explicitly describe a set of generators for the algebra $\mathcal{D}(W)$. This result 
shows the interplay between the differential operators associated with individual components of $W$ and the off-diagonal operators linking them.

\begin{cor}\label{generators}
  The algebra $\mathcal D(W)$ is generated by the set 
  $$\left  \{  \delta_{w_j} E_{j,j}\, , \mathcal T_{w_i,w_j} E_{i,j}, E_{j,j} \, : \, 1\leq i\neq j\leq N\right \}.$$
\end{cor}

 Now, we compute the algebra of eigenvalues $\Lambda(W) \subset \operatorname{Mat}_N(\mathbb{C}[n])$ for direct sums of classical scalar weights.

\begin{thm}\label{eigenv-directsum}
Let $w_{1}, \ldots, w_{N}$ be classical weights of the same type. For $W = w_{1} \oplus \cdots \oplus w_{N}$, the eigenvalues algebra $\Lambda(W)$ is given by
\[
\Lambda(W) = \sum_{i,j=1}^{N} \Lambda_{n}(\mathcal{T}{w_i,w_j}) \mathbb{C}[\Lambda_{n}(\delta_{j})] \, E_{i,j} \subset \operatorname{Mat}N(\mathbb{C}[n]),
\]
where $\mathbb{C}[\Lambda{n}(\delta_{j})]$ is as defined in \eqref{AE-escalar}, and
     \begin{enumerate}
    \item [(i) ] $\Lambda_{n}(\mathcal T_{w_j,w_j})=1,$ and   $\Lambda_{n}(\mathcal T_{w_1,w_2})=0$ if $w_1$ is not a Darboux transformation of $w_2$.
\item [(ii) ]  For Laguerre weights $w_\alpha$,  $k\in \mathbb N$,  we have 
$$ \Lambda_{n}(\mathcal T_{w_\alpha, w_{\alpha+k}}) = 1,  \qquad  \Lambda_{n}(\mathcal T_{w_\alpha+k, w_{\alpha}})= 
\frac{\Gamma(n+\alpha+k+1)}{\Gamma(n+\alpha+1)}.
$$
\item [(iii) ] For Jacobi weights $w_{\alpha,\beta}$, $k\in \mathbb{N}$, we have 
$$\Lambda_{n}(\mathcal T_{w_{\alpha,\beta+k}, w_{\alpha+k,\beta}}) = \frac{\Gamma(n+\beta+k+1)}{\Gamma(n+\beta+1)},  \qquad  \Lambda_{n}(\mathcal T_{w_{\alpha+k,\beta}, w_{\alpha,\beta+k}})=  \frac{\Gamma(n+\alpha+k+1)}{\Gamma(n+\alpha+1)}.$$

\end{enumerate}
\end{thm}
\begin{proof}
   The $(i,j)$-th entry of a differential operator $D \in  \mathcal D(W)$  
is a scalar   
differential operator  $ D_{i,j} \in \mathcal D(W_i, W_j)$. We have 
$$ \Lambda_n(D)=\begin{pmatrix}
    \Lambda_n(D_{1,1}) & \Lambda_n(D_{1,2}) & \cdots & \Lambda_n(D_{1,N})\\ \Lambda_n(D_{2,1}) & \Lambda_n(D_{2,2}) & \cdots & \Lambda_n(D_{2,N})
      \\  && \ddots  &&
\end{pmatrix}, $$
and $D_{i,j}= \mathcal{T}_{w_i,w_j} p_{i,j}(\delta_j)$, for some polynomial $p_{i,j}$. Then we have, 
$$\Lambda_n(D_{i,j})= \Lambda_n(\mathcal{T}_{w_i,w_j}) p_{i,j}(\Lambda_n(\delta_j)), $$ 
and the statement follows from Theorem \ref{T-operators}.
\end{proof}

To further understand the structure of the algebra $\mathcal{D}(W)$, where $W$ is a direct sum of classical scalar weights, we now focus on the computation of its center. 

We begin by analyzing the center of the algebra $\mathcal{D}(W)$ in the particular case where all scalar weights $w_j$ are Darboux equivalent. We consider weights of the form
$$W = w_{1} \oplus w_{2} \oplus \cdots \oplus w_{N}, \qquad \text{where } w_{i} \text{ is Darboux equivalent to } w_{j} \text{ for all } i, j.$$

\begin{thm}\label{center1}
Let $w_1, \dots, w_N$ be classical weights in the same Darboux equivalence class, and let $\delta_j$ denote the second-order differential operator associated with $w_j$, as given in Table \ref{weights-table}. Then, the center of the algebra $\mathcal{D}(w_{1} \oplus \cdots \oplus w_{N})$ is the polynomial algebra generated by the diagonal differential operator
\[
\Delta = \begin{psmallmatrix} \delta_{1} & & \\ & \ddots & \\ & & \delta_{N} \end{psmallmatrix}.
\]
Explicitly, we have
\[
\mathcal{Z}(w_{1} \oplus w_2 \oplus \cdots \oplus w_{N}) = \Big\{ p(\Delta) : p \in \mathbb{C}[x] \Big\}
= \left\{ \begin{pmatrix}
    p(\delta_1) & & \\ 
    & \ddots & \\
    & & p(\delta_N)
\end{pmatrix} : p \in \mathbb{C}[x] \right\}.
\]
\end{thm}
\begin{proof}
If $w_i$ is a Darboux transformation of $w_j$, then $\Lambda_n(\delta_i) = \Lambda_n(\delta_j)$, and consequently, we have $\Lambda_n(\Delta) = \Lambda_n(\delta_1) I$. Thus, for any polynomial $p \in \mathbb{C}[x]$, the differential operator $p(\Delta)$ belongs to the center $\mathcal{Z}(W)$.

Now let $D \in \mathcal{Z}(w_{1} \oplus w_{2} \oplus \cdots \oplus w_{N})$. Since $D$ commutes with $E_{i,i}$ for all $1 \leq i \leq N$, it follows that $D$ must be a diagonal matrix, and each diagonal entry $D_{i,i}$ belongs to $\mathcal{D}(w_i)$.

If $w_i$ is a Darboux transformation of $w_j$, there exists a nonzero scalar operator $\tau$ such that $\tau E_{i,j} \in \mathcal{D}(w_{1} \oplus w_{2} \oplus \cdots \oplus w_{N})$. This implies that $\Lambda_n(D)$ commutes with $\Lambda_n(\tau E_{i,j}) = \Lambda_n(\tau) E_{i,j}$ for all $n$. Consequently, we obtain $\Lambda_n(D_{i,i}) = \Lambda_n(D_{j,j})$ for all $n$.

Since $D_{i,i}$ is of the form $p_i(\delta_i)$ for some polynomial $p_i \in \mathbb{C}[x]$, the equality $\Lambda_n(D_{i,i}) = \Lambda_n(D_{j,j})$ implies that $p_i(x) = p_j(x)$ for all $i, j$. Therefore, $D = p(\Delta)$ for some $p \in \mathbb{C}[x]$, and the result is proved.
\end{proof}

\begin{example}
We consider the direct sum of two Laguerre weights
\[
W(x) = w_{\alpha}(x) \oplus w_{\alpha+1}(x),
\]
where the second-order differential operator  associated with $w_\alpha$ is given by
$\delta_{\alpha} = \partial^2 x + \partial (\alpha + 1 - x)$.

Since the weights \( w_{\alpha} \) and \( w_{\alpha+1} \) are Darboux equivalent, the center of the algebra \( \mathcal{D}(W) \) is generated by the diagonal operator
\[
\Delta  
= \begin{pmatrix} \partial^2 x + \partial (\alpha + 1 - x)  & 0 \\ 0 &\partial^2 x + \partial (\alpha + 2 - x)\end{pmatrix} .
\]
and the center is given by
$$
\mathcal{Z}(W) = \Big\{ p(\Delta) : p \in \mathbb{C}[x] \Big\}.$$
\end{example}

We now extend our analysis to the general case, where the scalar weights  $w_1, \dots, w_N$  may or may not be Darboux equivalent. While the final structure of the center remains similar, the key difference lies in the freedom to take distinct polynomials  $p_i(\delta_i)$  for each weight  $w_i$  when they are not Darboux equivalent.

By reordering the weights, we can assume that those which are Darboux equivalent are grouped together. This reordering does not alter the algebraic structure of the problem but allows us to present the result in a cleaner form, where the center appears as a block-diagonal structure, with each block corresponding to a group of Darboux-equivalent weights.

This result follows immediately from the case where all weights are Darboux equivalent and from Theorem \ref{algdirect}, which describes the algebra of differential operators for direct sums of classical weights.

\begin{thm}\label{cent direct}
Let \( W = W_{1} \oplus \cdots \oplus W_{m} \), where each \( W_{i} \) is a direct sum of classical scalar weights within the same Darboux equivalence class, i.e.,
\[
W_{i} = w_{i,1} \oplus \cdots \oplus w_{i,k_{i}},
\]
where \( w_{i,j} \sim w_{i,\ell} \) for all \( j, \ell \), and \( w_{i,\ell} \not\sim w_{r,j} \) for all \( \ell, j \) and \( i \neq r \). Let \( \delta_{i,j} \) denote the second-order differential operator associated with the weight \( w_{i,j} \). Then, the center of \( \mathcal{D}(W) \) is given by
\[
\mathcal{Z}(W) = \bigoplus_{i=1}^{m} \mathcal{Z}(W_{i}) = \mathbb{C}[\Delta_{1}] \oplus \cdots \oplus \mathbb{C}[\Delta_{m}],
\]
where each \( \Delta_{i} \) is the diagonal matrix
\[
\Delta_{i} = \begin{psmallmatrix} \delta_{i,1} & & \\ & \ddots & \\ & & \delta_{i,k_{i}} \end{psmallmatrix}.
\]
\end{thm}

We illustrate Theorem \ref{cent direct} with the direct sum of four Laguerre weights:
\[
W(x) = w_{\alpha-\frac{1}{2}}(x) \oplus w_{\alpha}(x) \oplus w_{\alpha+1}(x) \oplus w_{\alpha+\frac{1}{2}}(x), \qquad (\alpha > -\tfrac{1}{2}).
\]

To simplify the structure, we note that the weights naturally split into two groups of Darboux-equivalent weights
\begin{itemize}
    \item \( W_1(x) = w_{\alpha}(x) \oplus w_{\alpha+1}(x) \), where \( w_{\alpha} \sim w_{\alpha+1} \),
    \item \( W_2(x) = w_{\alpha+\frac{1}{2}}(x) \oplus w_{\alpha-\frac{1}{2}}(x) \), where \( w_{\alpha+\frac{1}{2}} \sim w_{\alpha-\frac{1}{2}} \).
\end{itemize}
However, the weights in \( W_1 \) are \textbf{not} Darboux equivalent to those in \( W_2 \).

Given this decomposition, there exists a nonsingular constant matrix \( M \) such that
\[
W(x) = M \, \Big( W_1(x) \oplus W_2(x) \Big) \, M^*.
\]

The centers of the algebras \( \mathcal{D}(W_1) \) and \( \mathcal{D}(W_2) \) are generated by the diagonal differential operators
\[
\Delta_1 = \begin{pmatrix} \delta_{\alpha} & 0 \\ 0 & \delta_{\alpha+1} \end{pmatrix}, \qquad
\Delta_2 = \begin{psmallmatrix} \delta_{\alpha+\frac{1}{2}} & 0 \\ 0 & \delta_{\alpha-\frac{1}{2}} \end{psmallmatrix},
\]
where \( \delta_{\alpha} \) is the second-order differential operator associated with the Laguerre weight \( w_{\alpha} \).
Thus, the center of the algebra \( \mathcal{D}(W_1 \oplus W_2) \) is
\[
\mathcal{Z}(W_1 \oplus W_2) = \mathcal{Z}(W_1 ) \oplus \mathcal Z(W_2)=  \left\{ 
\begin{pmatrix} p(\Delta_1) & 0 \\ 0 & q(\Delta_2) \end{pmatrix} \, : \, p, q \in \mathbb{C}[x] 
\right\}.
\]

Finally, since \( W \) is equivalent to \( W_1 \oplus W_2 \) up to the transformation \( M \), the center of \( \mathcal{D}(W) \) is
\[
\mathcal{Z}(W) = M \Big( \mathcal{Z}(W_1 \oplus W_2) \Big) M^{-1}.
\]

 \smallskip

We close this section by studying the structure of the algebra  $\mathcal{D}(W)$  in the case where the weight matrix  $W$  is a direct sum of Darboux-equivalent Laguerre weights.  We consider weight matrices of the form
$$W = w_{\alpha_1} \oplus w_{\alpha_2} \oplus \cdots \oplus w_{\alpha_N},$$
where the parameters  $\alpha_j$  are given by
\[   \alpha_j = \alpha + k_j , \quad \text { with } \; k_j \in \mathbb{Z} \quad  \text { and } \quad 0 = k_1 \leq k_2 \leq \cdots \leq k_N .\]

From Theorem \ref{algdirect}, the algebra \( \mathcal{D}(W) \) of differential operators can be written as
\[
\mathcal{D}(W) = \sum_{i,j=1}^{N} \mathcal{T}_{i,j} \cdot \mathcal{D}(w_{j}) \, E_{i,j},
\]
where the differential operators \( \mathcal{T}_{i,j} \) 
are given explicitly as
\begin{align*}
& \mathcal{T}_{i,j}   = 
(-\partial + 1)^{k_j - k_i}, 
& \mathcal{T}_{j,i}  = (\partial x + \alpha + k_j) \cdots (\partial x + \alpha + k_i + 1),
\end{align*}
for $k_i<k_j$ 
If  $k_i = k_j$, in particular when $i=j$ , both  $\mathcal{T}_{i,j}$  and  $\mathcal{T}_{j,i}$  are interpreted as the identity operator  1.
Let us observe that 
the degree of the operators \( \mathcal{T}_{i,j} \) and \( \mathcal{T}_{j,i} \) is \( k_j - k_i \).

From Corollary \ref{generators}, the algebra 
\( \mathcal{D}(W) \) is generated by the following set of elements
\[
D_j := \delta_{w_j} E_{j,j}, \quad F_{i,j} := \mathcal{T}_{i,j} E_{i,j}, \quad \text{and} \quad E_j := E_{j,j}, \quad \text{for } 1 \leq i \neq j \leq N.
\]
Here
\begin{itemize}
    \item \( D_j \) are the diagonal  
    differential operators associated with Laguerre  weights \( w_j \),
    \item \( E_j \) are diagonal idempotents satisfying
    \[
    I = E_1 \oplus \cdots \oplus E_N, \quad E_j^2 = E_j, \quad \text{and} \quad E_i E_j = 0 \quad \text{for } i \neq j,
    \]
    \item \( F_{i,j} \) are specific off-diagonal elements of the form \( \mathcal{T}_{i,j} E_{i,j} \).
\end{itemize}

It is important to note that the off-diagonal entries of the algebra \( \mathcal{D}(W) \) generally include 
elements of the form \( \mathcal{T}_{i,j} \cdot p(\delta_j) E_{i,j} \), where \( p \in \mathbb{C}[x] \). 
The operators $F_{i,j}$ satisfy 
$$F_{i,r}F_{s,j}=0, \text { for } i\neq j \quad \text{ and} \quad F_{i,r}F_{r,j}=q(\delta_i) F_{i,j}, \text{ for some } q\in \mathbb C[x]. $$
Moreover, 
$$ E_r F_{i,j}=0 \text { for } i\neq r \quad \text{ and} \quad E_i F_{i,j}=F_{i,j}= F_{i,j}E_j.$$

We can refine this set of generators and show that \( \mathcal{D}(W) \) can be generated by the smaller set
\[
\left\{ \Delta, F_{i,j} + F_{j,i}, E_j : 1 \leq i < j \leq N \right\},
\]
where \( \Delta \) is the sum of the diagonal generators
\[
\Delta := D_1 + D_2 + \cdots + D_N.
\]
This result follows directly from the relations satisfied by the operators \( F_{i,j} \).

More interestingly, the algebra \( \mathcal{D}(W) \) can be viewed as a finitely generated module over its center \( \mathcal{Z}(W) \). A generating set for this module is given by
\[
\left\{ F_{i,j}, E_j : 1 \leq i \neq j \leq N \right\}.
\]
This result highlights the algebraic structure of \( \mathcal{D}(W) \), where the off-diagonal elements \( F_{i,j} \) and the diagonal idempotents \( E_j \) provide a minimal generating set over the center.

\section{A family of irreducible  weight matrices: a Jacobi example}\label{sect-Jacobi}
In this section, we study an irreducible Jacobi-type weight matrix \( W \), which admits a second-order differential operator in its algebra \( \mathcal{D}(W) \). This example is significantly more challenging than the diagonal cases considered in the previous section due to the irreducibility of \( W \). 
To address these difficulties, we use the fact that \( W \) is Darboux-equivalent to a diagonal weight. This allows us to extract key information about the structure of \( \mathcal{D}(W) \), its center \( \mathcal{Z}(W) \), and its eigenvalue algebra \( \Lambda(W) \).

This weight matrix appeared previously in \cite{CGPSZ19}, as a solution of the matrix Bochner problem. After a suitable change of variables, it can be expressed as
\[
W(x) = (1-x)^{\alpha} (1+x)^{\beta} \left( W_2 \frac{(1-x)^2}{4} + W_1 \frac{(1-x)}{2} + W_0 \right), \qquad x \in (-1, 1),
\]
where the matrix coefficients \( W_2, W_1, W_0 \) are given by
\[
W_2 = \begin{psmallmatrix} \frac{v (v+2+\alpha+\beta)}{v + \alpha - \beta} & 0 \\ 0 & \frac{v(-v+2+\alpha+\beta)}{v-\alpha+\beta} \end{psmallmatrix}, \quad
W_1 = \begin{psmallmatrix} -(v + \alpha + \beta + 2) & (\alpha + \beta + 2) \\ (\alpha + \beta + 2) & -(-v + \alpha + \beta + 2) \end{psmallmatrix}, \quad
W_0 = \begin{psmallmatrix} \alpha + 1 & -\alpha - 1 \\ -\alpha - 1 & \alpha + 1 \end{psmallmatrix}.
\]
Here, the parameter $v$ satisfies \( |\alpha - \beta| < |v| < \alpha + \beta + 2 \). 

\smallskip

\begin{thm}\label{darb jac}
The weight \( W(x) \) is a Darboux transformation of the  diagonal 
weight \( \tilde{W}(x) = (1-x)^{\alpha+1}(1+x)^{\beta+1} I \).
\end{thm}

\begin{proof}
We consider the second-order differential operator \( D \in \mathcal{D}(\tilde{W}) \) given by
\begin{equation}\label{opertD-Jacobi}
    D = \begin{psmallmatrix} \delta_{\alpha+1,\beta+1} - \frac{1}{4}(-v+4+\alpha+\beta)(v+2+\alpha+\beta) & 0 \\ 
0 & \delta_{\alpha+1,\beta+1} - \frac{1}{4}(\alpha+\beta+v+4)(-v+2+\alpha+\beta) \end{psmallmatrix},
\end{equation}
where \( \delta_{\alpha+1,\beta+1} \) is the second-order differential operator associated with the Jacobi weight of parameters \( \alpha+1 \) and \( \beta+1 \), defined as
\[
\delta_{\alpha+1,\beta+1} = \partial^{2}(1-x^{2}) + \partial \big(\beta - \alpha - x(\alpha + \beta + 3)\big).
\]

We factorize \( D \) as \( D = \mathcal{V} \mathcal{N} \), where
\begin{equation}  \label{Voperator}
\mathcal{V} = \partial \begin{psmallmatrix} x + \frac{\beta-\alpha}{v} & 1 + \frac{\beta - \alpha}{v} \\ 1 + \frac{\alpha - \beta}{v} & x + \frac{\alpha - \beta}{v} \end{psmallmatrix} 
+ \begin{psmallmatrix} \frac{v+2+\alpha+\beta}{2} & 0 \\ 0 & \frac{-v+2+\alpha+\beta}{2} \end{psmallmatrix},
\end{equation}
and
\[
\mathcal{N} = \partial \begin{psmallmatrix} -x + \frac{\beta-\alpha}{v} & 1 + \frac{\beta-\alpha}{v} \\ 1 + \frac{\alpha-\beta}{v} & -x + \frac{\alpha-\beta}{v} \end{psmallmatrix} 
+ \begin{psmallmatrix} -\frac{1}{2}(-v+4+\alpha+\beta) & 0 \\ 0 & -\frac{1}{2}(\alpha+\beta+v+4) \end{psmallmatrix}.
\]

The operators \( \mathcal{V} \) and \( \mathcal{N} \) satisfy the relation
\[
\mathcal{N} = W \mathcal{V}^{\ast} \tilde{W}^{-1},
\]
and 
\[
\langle P \cdot \mathcal{V}, Q \rangle_{W} = \langle P, Q \cdot \mathcal{N} \rangle_{\tilde{W}} \quad \text{for all } P, Q \in \operatorname{Mat}_{2}(\mathbb{C}[x]).
\]
By Proposition \ref{strong1}, it follows that \( W \) is a Darboux transformation of \( \tilde{W} \). 
\end{proof}

From the factorization established in the previous theorem, we now derive an explicit sequence of orthogonal polynomials for \( W \).

\begin{cor}
A sequence of orthogonal polynomials with respect to the weight \( W \) can be expressed as  
\begin{equation*}
Q_n(x) = J_n^{(\alpha+1, \beta+1)}(x) \cdot \mathcal{V},
\end{equation*}
where \( \mathcal{V} \) is the degree-preserving operator defined in \eqref{Voperator}, and \( J_n^{(\alpha+1, \beta+1)}(x) \) denotes the sequence of monic orthogonal polynomials associated with the Jacobi weight \( w_{\alpha+1, \beta+1}(x) = (1-x)^{\alpha+1} (1+x)^{\beta+1} \).  
\end{cor}

Since \( W \) is a Darboux transformation of the diagonal weight \( \tilde{W} \), we can describe the algebra \( \mathcal{D}(W) \) explicitly. Using the Darboux framework and the results developed in this paper, we determine a set of generators for \( \mathcal{D}(W) \).

\begin{thm}\label{alg jacobi W}
    The algebra $\mathcal{D}(W)$ associated with the irreducible Jacobi-type weight  $W$ is generated by the set  $\{I, D_{1}, D_{2}, D_{3}, D_{4} \}$, where  the differential operators  $D_i$  are explicitly given by
    \begin{align*}
            D_{1} = \mathcal{N}\begin{psmallmatrix} 1 & 0 \\0 & 0 \end{psmallmatrix} \mathcal{V} & = \partial^{2} \begin{psmallmatrix}-x^{2} + \frac{(\alpha-\beta)^{2}}{v^{2}} & \frac{(vx+\alpha-\beta)(-v+\alpha-\beta)}{v^{2}} \\ \frac{(vx-\alpha+\beta)(v-\beta+\alpha)}{v^{2}} & 1 - \frac{(\alpha-\beta)^{2}}{v^{2}}\end{psmallmatrix} \\
            & \quad + \partial \begin{psmallmatrix} -x(\alpha+\beta+4) - \frac{(v-2)(\alpha-\beta)}{2v} & \frac{(-v+\alpha-\beta)(\alpha\beta-v+6)}{2v} \\ \frac{(v-\beta+\alpha)(v+2+\alpha+\beta)}{2v} & 0 \end{psmallmatrix} + \begin{psmallmatrix} -\frac{(-v+4+\alpha+\beta)(v+2+\alpha+\beta)}{4} & 0 \\ 0 & 0  \end{psmallmatrix}, \displaybreak[0]\\ 
            D_{2} = \mathcal{N}\begin{psmallmatrix} 0 & 1 \\0 & 0 \end{psmallmatrix} \mathcal{V}  & = \partial^{2} \begin{psmallmatrix}\frac{(vx + \alpha - \beta)(-v+\beta-\alpha)}{v^{2}} & - \frac{(vx + \alpha - \beta)^{2}}{v^{2}} \\ \frac{(v-\beta+\alpha)^{2}}{v^{2}} & \frac{(vx+\alpha-\beta)(v-\beta+\alpha)}{v^{2}} \end{psmallmatrix} \\
            & \quad + \partial \begin{psmallmatrix} \frac{(-v+\beta-\alpha)(\alpha+\beta-v+6)}{2v} & - \frac{(-v+4+\alpha+\beta)(vx+\alpha-\beta)}{v} \\ 0 & \frac{(v-\beta+\alpha)(-v+2+\alpha+\beta)}{2v} \end{psmallmatrix} + \begin{psmallmatrix} 0 & - \frac{(-v+4+\alpha+\beta)(-v+2+\alpha+\beta)}{4}\end{psmallmatrix}, \displaybreak[0]\\
            D_{3} = \mathcal{N}\begin{psmallmatrix} 0 & 0 \\1 & 0 \end{psmallmatrix} \mathcal{V} & = \partial^{2} \begin{psmallmatrix} \frac{(-vx+\alpha-\beta)(-v+\alpha-\beta)}{v^{2}} & \frac{(-v+\alpha-\beta)^{2}}{v^{2}} \\ - \frac{(-vx + \alpha - \beta)}{v^{2}} & \frac{(-vx+\alpha-\beta)(v-\alpha+\beta)}{v^{2}}\end{psmallmatrix} \\ 
            & \quad + \partial \begin{psmallmatrix} - \frac{(v+2+\alpha+\beta)(-v+\alpha-\beta)}{2v} & 0 \\ \frac{(\alpha + \beta + v + 4)(-vx + \alpha - \beta)}{v} & \frac{(\alpha + \beta + v + 6)(-v + \alpha - \beta)}{2v}\end{psmallmatrix} + \begin{psmallmatrix}0 & 0 \\ -\frac{(\alpha+\eta+v+4)(v+2+\alpha+\beta)}{4} & 0 \end{psmallmatrix}, \displaybreak[0]\\
             D_{4} = \mathcal{N}\begin{psmallmatrix} 0 & 0 \\0 & 1 \end{psmallmatrix} \mathcal{V} & = \partial^{2} \begin{psmallmatrix} 1 - \frac{(\alpha - \beta)^{2}}{v^{2}} & \frac{(vx + \alpha - \beta)(v + \beta - \alpha)}{v^{2}} \\ \frac{(-vx+\alpha-\beta)(v-\beta+\alpha)}{v^{2}} & -x^{2} + \frac{(\alpha-\beta)^{2}}{v^{2}}\end{psmallmatrix} \\
             & \quad + \partial \begin{psmallmatrix} 0 & \frac{(v - \alpha + \beta)(-v + 2 + \alpha + \beta)}{2v} \\ - \frac{(v - \beta + \alpha)(\alpha + \beta + v + 6)}{2v} & -x(\alpha + \beta + 4)- \frac{(v+2)(\alpha-\beta)}{v}\end{psmallmatrix} + \begin{psmallmatrix} 0 & 0 \\ 0 & - \frac{(\alpha+\beta+v+4)(-v+2+\alpha+\beta)}{4}\end{psmallmatrix}.
    \end{align*}
\end{thm}

\begin{proof}
The algebra \( \mathcal{D}(\tilde{W}) \) contains only even-order differential operators, and since both \( \mathcal{V} \) and \( \mathcal{N} \) have nonsingular leading coefficients, it follows that all differential operators in \( \mathcal{D}(W) \) must also be of even order. It is straightforward to verify that the only operators of order zero are of the form \( cI \), where \( c \in \mathbb{C} \).

Define the operator \( D_{5} = D_{1} + D_{4} = \mathcal{N}\mathcal{V} \). Explicitly, we have
\begin{equation*}
    \begin{split}
        D_{5} &= \partial^{2}(1-x^{2})I + \partial \begin{psmallmatrix} -x(\alpha+\beta+4)- \frac{(v-2)(\alpha-\beta)}{v} & \frac{2(-v+\alpha-\beta)}{v} \\ - \frac{2(v-\beta+\alpha)}{v} & -x(\alpha+\beta+4)- \frac{(v+2)(\alpha-\beta)}{v} \end{psmallmatrix} \\ 
        &\quad + \begin{psmallmatrix} - \frac{(-v+4+\alpha+\beta)(v+2+\alpha+\beta)}{4} & 0 \\ 0 & \frac{(\alpha+\beta+v+4)(-v+2+\alpha+\beta)}{4} \end{psmallmatrix}.
    \end{split}
\end{equation*}

Now, let \( B = \sum_{j=0}^{2m} \partial^{j} F_{j}(x) \in \mathcal{D}(W) \), with \( m \geq 1 \). We observe that the composition \( \mathcal{V} B \mathcal{N} \) has order \( 2m+2 \) and belongs to \( \mathcal{D}(\tilde{W}) \). Its leading coefficient is given by
\[
V_{1}(x) F_{2m}(x) N_{1}(x) = (1-x^{2})^{m+1} \begin{pmatrix} k_{1} & k_{2} \\ k_{3} & k_{4} \end{pmatrix},
\]
for some \( k_{1}, k_{2}, k_{3}, k_{4} \in \mathbb{C} \), where \( V_{1}(x) \) and \( N_{1}(x) \) are the leading coefficients of \( \mathcal{V} \) and \( \mathcal{N} \), respectively
\begin{equation}\label{lead darb}
    V_{1}(x) = \begin{psmallmatrix} x + \frac{\beta-\alpha}{v} & 1 + \frac{\beta - \alpha}{v} \\ 1 + \frac{\alpha - \beta}{v} & x + \frac{\alpha - \beta}{v} \end{psmallmatrix}, \quad 
    N_{1}(x) = \begin{psmallmatrix} -x + \frac{\beta-\alpha}{v} & 1 + \frac{\beta-\alpha}{v} \\ 1 + \frac{\alpha-\beta}{v} & -x + \frac{\alpha-\beta}{v} \end{psmallmatrix}.
\end{equation}

From the above, the leading coefficient \( F_{2m}(x) \) of \( B \) is given by
\begin{equation} \label{lead coeff}
    F_{2m}(x) = V_{1}(x)^{-1} (1-x^{2})^{m+1} \begin{pmatrix} k_{1} & k_{2} \\ k_{3} & k_{4} \end{pmatrix} N_{1}(x)^{-1}.
\end{equation}

Now, consider the differential operator
\[
E = D_{5}^{m-1} \left( k_{1}D_{1} + k_{2}D_{2} + k_{3}D_{3} + k_{4}D_{4} \right).
\]
The operator \( E \) has the same leading coefficient as \( B \), and thus \( B - E \in \mathcal{D}(W) \) is a differential operator of order at most \( 2m-2 \). 

By applying this process inductively, we conclude that \( B \) is generated by \( \{ I, D_1, D_2, D_3, D_4 \} \), completing the proof.
\end{proof}

To further analyze the algebra \( \mathcal{D}(W) \) and establish the relations between the differential operators \( D_{1}, D_{2}, D_{3}, D_{4} \), we now compute the eigenvalue maps \( \Lambda_n(D_i) \) for \( i = 1, 2, 3, 4 \). These explicit expressions will play a central role in proving the upcoming proposition.

The eigenvalue matrices \( \Lambda_n(D_i) \), which encode the action of the operators \( D_i \) on the monic orthogonal polynomials associated with \( W \), are given by
\[
\begin{aligned}
    \Lambda_{n}(D_{1}) &= \begin{psmallmatrix} -\frac{1}{4}(2n-v+4+\alpha+\beta)(2n+v+2+\alpha+\beta) & 0 \\ 0 & 0 \end{psmallmatrix}, \\
    \Lambda_{n}(D_{2}) &= \begin{psmallmatrix} 0 & -\frac{1}{4}(2n-v+4+\alpha+\beta)(2n-v+2+\alpha+\beta) \\ 0 & 0 \end{psmallmatrix}, \\
    \Lambda_{n}(D_{3}) &= \begin{psmallmatrix} 0 & 0 \\ -\frac{1}{4}(2n+v+4+\alpha+\beta)(2n+v+2+\alpha+\beta) & 0 \end{psmallmatrix}, \\
    \Lambda_{n}(D_{4}) &= \begin{psmallmatrix} 0 & 0 \\ 0 & -\frac{1}{4}(2n+v+4+\alpha+\beta)(2n-v+2+\alpha+\beta) \end{psmallmatrix}.
\end{aligned}
\]

Using these expressions, we can deduce the algebraic relations satisfied by the generators  $D_{1}$, $D_{2},$ $ D_{3},$ $ D_{4}$  in  $\mathcal{D}(W)$, which are presented in the following proposition.

\begin{prop}\label{relations}
    The differential operators  $D_{1},D_{2},D_{3},D_{4}$, introduced in Theorem \ref{alg jacobi W}, satisfy the following relations
   $$ \begin{aligned}
    & D_{2}D_{3} = D_{1}^{2} + vD_{1}, \quad D_{3}D_{2} = D_{4}^{2} - vD_{4}, \quad 
      D_{1}D_{2} - D_{2}D_{4} = -vD_{2}, \quad D_{4}D_{3} - D_{3}D_{1} = vD_{3}, \\
       & D_{1}D_{4} = D_{4}D_{1}= 0, \quad D_{1}D_{3} = 0, 
    \quad D_{4}D_{2} = 0,  \quad 
    D_{2}D_{1} = 0, \quad D_{2}^{2} = 0, \quad D_{3}D_{4} = 0, \quad D_{3}^{2} = 0.
\end{aligned}
$$
\end{prop}

\smallskip

To determine the center \( \mathcal{Z}(W) \) of the algebra \( \mathcal{D}(W) \), we proceed by leveraging the results of Theorem \ref{darb center}. Recall that this theorem guarantees the existence of a differential operator in the center of \( \mathcal{D}(W) \), which can be factorized as a product \( \mathcal{V}\mathcal{N} \) under appropriate conditions. This factorization plays a fundamental role in connecting the centers of \( \mathcal{D}(W) \) and \( \mathcal{D}(\tilde{W}) \), as it establishes the inclusion
\[
\mathcal{N} \mathcal{Z}(\tilde{W}) \mathcal{A} \mathcal{V} \subset \mathcal{Z}(W).
\]
This inclusion allows us to exploit the knowledge of the center \( \mathcal{Z}(\tilde{W}) \), which is already known since \( \tilde{W} \) is a diagonal Jacobi weight. Specifically, Theorem \ref{center1} asserts that \( \mathcal{Z}(\tilde{W}) \) is generated by the second-order differential operator \( \tilde{D} = \delta_{\alpha+1,\beta+1} I \), the classical Jacobi operator.

Following this approach, we consider the fourth-order differential operator \( \mathcal{D} \in \mathcal{Z}(\tilde{W}) \), defined as
\[
\mathcal{D} = \mathcal{A} D, \quad \text{where} \quad 
\mathcal{A} = \begin{psmallmatrix} \delta_{\alpha+1,\beta+1} - \frac{1}{4}(\alpha+\beta+v+4)(-v+2+\alpha+\beta) & 0 \\ 0 & \delta_{\alpha+1,\beta+1} - \frac{1}{4}(-v+4+\alpha+\beta)(v+2+\alpha+\beta)  \end{psmallmatrix}.
\]
The eigenvalue map of \( \mathcal{D} \) satisfies $\Lambda_{n}(\mathcal{D}) = p(n) I$,  where 
\[
p(n) = \frac{1}{16} (2n -v +2 +\alpha + \beta)(2n + \alpha + \beta +v + 4)(2n + v +2 +\alpha +\beta)(2n-v+4+\alpha+\beta).
\]
Thus, \( \mathcal{D} \) belongs to \( \mathcal{Z}(\tilde{W}) \) and can be factorized as \( \mathcal{D} = (\mathcal{A}\mathcal{V})\mathcal{N} \).

By combining these results with Proposition \ref{cent}, we find that the fourth-order and sixth-order differential operators
\[
Z_{1} = \mathcal{N} \mathcal{A} \mathcal{V} \quad \text{and} \quad Z_{2} = \mathcal{N} \delta_{\alpha+1,\beta+1} \mathcal{A} \mathcal{V}
\]
belong to the center \( \mathcal{Z}(W) \). These operators can be expressed in terms of the generators \( D_{1}, D_{4} \) of \( \mathcal{D}(W) \), as follows
\begin{equation}\label{genZW}
    \begin{split}
  Z_{1} &= D_{1}^{2} + D_{4}^{2} + v(D_{1} - D_{4}), \\
    Z_{2} &= D_{1}^{3} + D_{4}^{3} + \frac{1}{4}(-v+4+\alpha+\beta)(v+2+\alpha+\beta)(D_{1}^{2} + vD_{1}) \\
    &\quad + \frac{1}{4}(-v+2+\alpha+\beta)(v+4+\alpha+\beta)(D_{4}^{2} - vD_{4}) + v(D_{1}^{2} - D_{4}^{2}).
    \end{split}
\end{equation}

The leading coefficients of \( Z_{1} \) and \( Z_{2} \) are \( (1-x^2)^{2}I \) and \( (1-x^2)^{3}I \), respectively.

\smallskip

We are now ready to describe the center \( \mathcal{Z}(W) \) of the algebra \( \mathcal{D}(W) \). Using the structure provided by the Darboux transformation and the explicit expressions for the central operators \( Z_{1} \) and \( Z_{2} \), we can state the following result.

\begin{thm}
The center \( \mathcal{Z}(W) \) of the algebra \( \mathcal{D}(W) \) is generated by \( \{ I, Z_{1}, Z_{2} \} \), where \( Z_{1} \) and \( Z_{2} \) are the fourth-order and sixth-order differential operators, respectively, as given in \eqref{genZW}. Explicitly, we have
\[
\mathcal{Z}(W) = \mathbb{C}I \oplus \mathcal{N} \mathcal{Z}(\tilde{W}) \mathcal{A} \mathcal{V}.
\]
\end{thm}

\begin{proof}
Let \( B \in \mathcal{Z}(W) \) be a differential operator of order \( 2m \). By Theorem \ref{alg jacobi W}, every operator in \( \mathcal{D}(W) \) can be written as
\[
B = k_{1}D_{1} + k_{2}D_{2} + k_{3}D_{3} + k_{4}D_{4} + rI, \quad \text{for some } k_{1}, k_{2}, k_{3}, k_{4}, r \in \mathbb{C}.
\]
To belong to the center \( \mathcal{Z}(W) \), \( B \) must commute with each generator \( D_{i} \) (for \( i = 1, \dots, 4 \)). This condition \( D_{i}B = BD_{i} \) implies that \( k_{1} = k_{2} = k_{3} = k_{4} = 0 \). Therefore, there are no second-order differential operators in \( \mathcal{Z}(W) \).

Now, let \( B \in \mathcal{Z}(W) \) be a higher-order operator, i.e., of order \( 2m \) with \( m > 1 \). By Proposition \ref{cent}, the operator \( \mathcal{A} \mathcal{V} B \mathcal{N} \) belongs to \( \mathcal{Z}(\tilde{W}) \), where \( \mathcal{Z}(\tilde{W}) = \mathbb{C}[\delta_{\alpha+1, \beta+1}I] \). This operator has order \( 2m + 4 \) and its leading coefficient satisfies
\begin{align*}
F_{2m}(x) & = \frac{1}{(1-x^{2})} V_{1}(x)^{-1} 
\begin{pmatrix} a(1-x^{2})^{m+2} & 0 \\ 0 & a(1-x^{2})^{m+2} \end{pmatrix} N_{1}(x)^{-1} \\
& = \begin{pmatrix} a(1-x^{2})^{m} & 0 \\ 0 & a(1-x^{2})^{m} \end{pmatrix},
\end{align*}
for some \( a \in \mathbb{C} \).

Since the leading coefficient of \( B \) is proportional to \( (1-x^2)^{m} \), there exist non-negative integers \( s \) and \( t \) such that \( 4t + 6s = 2m \). Therefore, the differential operator \( B - aZ_{1}^{t}Z_{2}^{s} \) is of order at most \( 2m - 2 \) and remains in \( \mathcal{Z}(W) \).

By induction on the order \( 2m \), we conclude that \( B \) can be generated by \( I, Z_{1}, \) and \( Z_{2} \), completing the proof.
\end{proof}

To complete the description of the center \( \mathcal{Z}(W) \), we study the algebraic relations satisfied by its generators. 
The following proposition establishes the key relation between \( Z_{1} \) and \( Z_{2} \).

\begin{prop}
The generators \( Z_{1} \) and \( Z_{2} \) of the center \( \mathcal{Z}(W) \) satisfy the following algebraic relation
\begin{equation*}
    \begin{split}
         0 & = Z_{1}^{3}-Z_{2}^{2} + \left(\frac{(4+\alpha+\beta)(\alpha+\beta+2)-v^{2}}{2} \right )Z_{1}Z_{2} \\
        & \qquad - \frac{1}{16}(v+2+\alpha+\beta)(-v+2+\alpha+\beta)(v+4+\alpha+\beta)(-v+4+\beta+\alpha)Z_{1}^{2}.
    \end{split}
\end{equation*}
\end{prop}

To complete our study of the algebra \( \mathcal{D}(W) \), we now focus on its structure as a module over its center \( \mathcal{Z}(W) \). It is known from \cite{CY18} that the algebra \( \mathcal{D}(W) \) is always finitely generated as a module over its center. In the following theorem, we determine an explicit set of generators for \( \mathcal{D}(W) \), viewed as a \( \mathcal{Z}(W) \)-module.

\begin{thm}
    Let \( D_{1}, D_{2}, D_{3}, D_{4} \) be the differential operators defined in Theorem \ref{alg jacobi W}, and let \( D_{5} = D_{1} + D_{4} \). The algebra \( \mathcal{D}(W) \) is generated as a \( \mathcal{Z}(W) \)-module by the set 
    $$ \{ I, D_{j}, D_5 D_{j} : \, j=1,2,3,4 \}.$$
\end{thm}
\begin{proof}
    Let \( G_{i}(x) \) denote the leading coefficient of \( D_{i} \) for \( i = 1,2,3,4 \). By definition, the leading coefficient of \( D_{5}D_{i} \) is \( (1-x^2)G_{i}(x) \). Recall that the leading coefficients of the central operators \( Z_{1} \) and \( Z_{2} \) are \( (1-x^2)^2 \) and \( (1-x^2)^3 \), respectively.

    Let \( B = \sum_{j=0}^{2m} \partial^j F_{j}(x) \in \mathcal{D}(W) \) be a differential operator of order \( 2m \). By \eqref{lead coeff}, the leading coefficient \( F_{2m}(x) \) can be expressed as a linear combination of \( G_{1}(x), G_{2}(x), G_{3}(x), G_{4}(x) \), scaled by \( (1-x^2)^{m-1} \). We write
    \[
    F_{2m}(x) = (1-x^2)^{m-1} \left( k_{1}G_{1}(x) + k_{2}G_{2}(x) + k_{3}G_{3}(x) + k_{4}G_{4}(x) \right),
    \]
    for some \( k_{1}, k_{2}, k_{3}, k_{4} \in \mathbb{C} \).

    If \( m \leq 2 \), then \( B \) is a linear combination of the operators in \( E \), as the highest order terms match those in the set \( E \). Now, assume \( m > 2 \) and proceed by induction. There exist nonnegative integers \( s, t \) such that \( 2s + 3t = m-1 \). We construct the differential operator
    \[
    B' = Z_{1}^{s} Z_{2}^{t} \left( k_{1}D_{1} + k_{2}D_{2} + k_{3}D_{3} + k_{4}D_{4} \right),
    \]
    where \( Z_{1} \) and \( Z_{2} \) are the central elements of order 4 and 6, respectively.

    The operator \( B - B' \) has order at most \( 2m-2 \). By the inductive hypothesis, \( B - B' \) can be expressed as a linear combination of the generators in \( E \). Thus, the statement holds for all \( m \geq 0 \).
\end{proof}

\section*{Funding}
This work was partially supported by SeCyT-UNC and CONICET, Grant PIP 1220150100356.

\

\end{document}